\documentclass[leqno]{amsart}
\usepackage{amsmath,amssymb}

\usepackage[all]{xy}

\begin{document}
\begin{center}
\title{Non-linear evolution equations and hyperelliptic covers of elliptic curves}
\author{Armando Treibich}

\maketitle
\end{center}
  
\noindent \medskip\section{ Introduction}

\noindent\textbf{1.1.}
A huge variety of nonlinear integrable processes and phenomena in physics and mathematics can be described by a few nonlinear partial derivative equations (e.g.: \textit{Korteweg-deVries} and \textit{Kadomtsev-Petviashvili, 1D and 2D Toda, sine-Gordon, nonlinear Schr\"odinger}). For almost 40 years a full range of methods coming from distinct areas were developped in order to deal and present exact solutions of the latter equations (e.g.: \cite{A-McK-M} till \cite{Z-S} and their references). Zero-curvature equations, Lax pair's presentation and inverse scattering methods revolutionized the whole domain (\cite{Lax}, \cite{Z-S}). Rational and trigonometric exact solutions (\cite{A-McK-M}, \cite{C}, \cite{Hi}) were followed by quasi-periodic ones, also called \textit{finite-gap}, given in terms of the theta function of an arbitrary hyperelliptic curve, via the Its-Matveev formula or its variants (\cite{D-M-N}, \cite{I-M}). A few years later I.M. Krichever made a major contribution in  \cite{K.1}, extending the latter results to finite-gap solutions of the \textit{KP equation} associated to an arbitrary compact Riemann surface. M. Sato's infinite dimensional approach, developped in the beginning of the 80's (\cite{Sato}, \cite{S-W},\cite{K-M}), further generalized Krichever's dictionnary as well as the classical theta and Baker-Akhiezer function. From then on, all previously studied non-linear evolution equations were reconsidered, and considerable effort was made in order to find doubly periodic solutions to each one of them. The starting point to this new trend was Krichevers's seminal article \cite{K.2}. The first doubly periodic solutions to the \textit{KdV equation} and a remarkable connection with the elliptic \textit{Calogero-Moser integrable system} had already been found (e.g.: \cite{A-McK-M} \& \cite{D-N}, as well as \cite{C} for the rational/trigonometric case), but \cite{K.2} generalizes to an equivalence between the elliptic C-M integrable system and the KP solutions, doubly periodic in $x$. More precisely, given $n\geq 1$ and the lattice $L \subset \mathbb{C}$, the corresponding elliptic Calogero-Moser integrable system is solved. Its ($2n$-dimensional) phase space is cut out by the Jacobian Varieties of an $n$-dimensional family of genus $n$ marked compact Riemann surfaces, each one of which  (is effectively constructed and) gives rise to  KP solutions $L$-periodic in $x\in \mathbb{C}$. The analogous problems for the \textit{KdV}, \textit{1D Toda}, \textit{NL Schr\"odinger}, \textit{sine-Gordon} equations and related problems (\cite{McK-VM}, \cite{M}, \cite{P}) amount to finding hyperelliptic curves equipped with a projection onto $X$, satisfying specific geometrical properties, as briefly explained hereafter.\\
 Let indeed $\pi: (\Gamma,p) \to (X,q)$ be an arbitrary ramified cover, where $\pi(p)=q$ and $(X,q) $ is the elliptic curve $(\mathbb{C}/L,0)$. Up to a translation, there exist canonical copies of $\Gamma$ and $X$ inside $Jac\,\Gamma$, the Jacobian variety of $\Gamma$. Consider the flag $\{0\} \subsetneq V^1_{\Gamma,p}\ldots \subsetneq V^g_{\Gamma,p}$, of hyperosculating spaces to $\Gamma$ at $p$, and $T_o X $  the tangent line to (the copy of ) $X$, inside $Jac\,\Gamma$. \\

The $d$-th case of the \textbf{KP equation}: $\frac{3}{4}u_{yy}+\frac{\partial}{\partial x}\big(u_t+\frac{1}{4}(6uu_x\,$-$\,u_{xxx})\big)$.\\ 

We will call $\pi:(\Gamma,p) \to (X,q)$ a \textit{d-osculating cover} if $T_o X \subset V^d_{\Gamma,p}\backslash V^{d-1}_{\Gamma,p}$. Such covers, studied and constructed for any $d\geq 1$, give rise to KP solutions $L$-periodic with respect to the $d$-th KP flow (cf. \cite{T-V.1} for  $d=1$ and \cite{T.2} for any other $d$). \\

The $d$-th case of the \textbf{KdV equation}: $u_t+\frac{1}{4}(6uu_x\,$-$\,u_{xxx})$. \\

Recall that $p\in \Gamma$ is a \textit{Weierstrass point} of the \textit{hyperelliptic curve} $\Gamma$, if and only if there exists a degree-$2$ projection $\Gamma \to \mathbb{P}^1$, ramified at $p$. Or in other words, if and only if there exists  an involution, say $\tau_\Gamma :\Gamma \to \Gamma$,  fixing $p$ and such that the quotient curve $\Gamma /\tau_{\Gamma}$ is isomorphic to $\mathbb{P}^1$. Let $\pi: (\Gamma,p) \to (X,q)$ be a \textit{d-osculating cover} such that $\Gamma$ is \textit{hyperelliptic} and $p \in \Gamma$ a \textit{Weierstrass point}. Then, all KdV solutions  classically associated to  $(\Gamma,p)$ are $L$-periodic with respect to the $d$-th KdV flow.\\

The \textbf{Non Linear Schr\"odinger}: $\quad ip_y+p_{xx} \mp  8|p|^2p=0$\\

and the \textbf{1D Toda case}:  $\quad \quad\frac{\partial^2}{\partial t^2}\varphi_n=exp(\varphi_n\,$-$\,\varphi_{n\textnormal{-}1}) \,$-$\,exp(\varphi_{n+1}\,$-$\,\varphi_n) $.\\

Let $\pi: (\Gamma,p^+) \to (X,q)$ be a \textit{1-osculating cover} (i.e.: also called a \textit{tangential cover} in \cite{T.1}) such that $\Gamma$ is \textit{hyperelliptic} and $p^+ \in \Gamma$ is not a \textit{Weierstrass point}. Then, all nonlinear Schr\"odinger \& 1D Toda solutions classically associated to  $(\Gamma,p ^+, \tau_\Gamma(p^+))$, are $L$-periodic in $x$ and in $t$, respectively.\\

The \textbf{sine-Gordon case}: $\quad u_{xx}\,$-$\,u_{tt}=sin u$.\\

Let $\Gamma$ be a \textit{hyperelliptic} curve, equipped with a projection $\pi:\Gamma \to X$ and two \textit{Weierstrass points}, say $p,p'\in \Gamma$, such that the tangent line $T_oX$ is contained in the plane $V^1_{\Gamma,p}+V^1_{\Gamma,p'}$, generated by the tangents to $\Gamma$ at $p$ and $p'$ (inside $Jac\,\Gamma$). Then, up to choosing suitable local coordinates of $\Gamma$ at $p$ and $p'$, the sine-Gordon solutions classically associated to $(\Gamma,p,p')$ are $L$-periodic in $x$.\\

 The \textbf{KP} case being rather well understood, we will focus on the three other cases, and in particular, on ramified projections $\pi :\Gamma \rightarrow X$, of a hyperelliptic curve onto the fixed elliptic one, marked at, either one or two \textit{Weierstrass points} (\textbf{KdV} and \textbf{sine-Gordon} cases), or two points exchanged by the hyperelliptic involution. Studying the tangent and osculating spaces at the marked points (in $Jac\,\Gamma$) is an interesting geometric problem which, I believe, does not need any further motivation. It was first considered, however, through its links with $L$-periodic solutions of the \textit{Korteweg-deVries} equation (e.g.: \cite{A-McK-M}, \cite{D-N}, \cite{I-M}, \cite{K.2}, \cite{S.1}, \cite{T-V.1} for $d=1$ and \cite{S.3}, \cite{A-K-V}, \cite{F}, \cite{F-T} for $d = 2$), as well as the \textbf{Toda}, \textbf{sine-Gordon} and \textbf{nonlinear Schr\"odinger} equations (e.g.: \cite{S.2}, \cite{B}, \cite{S.4}). 
Studying their general properties (such as the relations between the genus and the degree of the cover), and constructing examples in any genus, will be the main issues of this article.

After fixing a lattice $L \subset \mathbb{C}$ defining the marked elliptic curve $(X,q):=(\mathbb{C}/L , 0)$, we will develop in section \textbf{3} a well suited algebraic-surface approach, for studying the structure of all ramified covers of $X$ we are interested in, and their canonical factorization through a particular algebraic surface. Natural numerical invariants will then be defined, in terms of which we will characterize the latter covers and, ultimately, construct arbitrarily high genus examples to each case. \medskip  \\
\noindent\textbf{1.2.}
 We sketch hereafter the structure and main results of our article.

  \begin {enumerate}
  
 \item We start section \textbf{2}defining the Abel rational embedding of a curve $\Gamma$, of positive genus $g$, into its \textit{generalized Jacobian}, $Jac\,\Gamma$,  and construct the \textit{flag of hyperosculating spaces} $\{0\} \subsetneq  V_{1,p}\ldots \subsetneq V_{g,p} = H^1(\Gamma,O_\Gamma)$, at the image of any smooth point $p \in \Gamma$. From then on, we restrict to Jacobians of hyperelliptic curves such that $Jac\,\Gamma$ contains the elliptic curve $(X,q)=(\mathbb{C}/L, 0)$, or equivalently, to any \textit{hyperelliptic cover} $\pi: (\Gamma,p) \to (X,q)$. Dualizing  such a cover $\pi$, we obtain a homomorphism $\iota_\pi:X \to Jac\, \Gamma$, with image an elliptic curve isogeneous to $X$. Let $d$ be the smallest positive integer, called the \textit{osculating order} of $\,\pi$, such that the tangent line defined by $\iota_\pi(X)$ is contained in $V_{d,p} \subset H^1(\Gamma,O_\Gamma)$. Whenever $p\in \Gamma$ is a \textit{Weierstrass} point, $\pi$ is called a \textit{hyperelliptic $d$-osculating cover}, and gives rise to  KdV solutions, $L$-periodic with respect to the $d$-th KdV flow. Such covers are characterized by the existence of a particular projection $\kappa:\Gamma \rightarrow \mathbb{P}^1$ (\textbf{2.6.}). Given any \textit{hyperelliptic cover} $\pi$, marked at, either two points exchanged by the hyperelliptic involution, or two \textit{Weierstrass} points, we also find analogous characterizations for $\pi$ to solve, the \textbf{NL Schr\"odinger \& 1D Toda} or the \textbf{sine-Gordon case} (\textbf{2.9., 2.10.}).\medskip 
 
\item The latter characterizations \textbf{2.6.} pave the way to the algebraic surface approach
  developed in the remaining sections. The main characters are played by
   three projective surfaces and corresponding morphisms, canonically associated to $X$:
    
   \medskip $\centerdot \,\,\pi_S : S \rightarrow X$ :\,\,a particular ruled surface;
   
   \medskip $\centerdot \,\,e:S^\bot \to S$ : \,\,the blow-up of $S$, at the $8$ fixed points of its involution; 
   
   \medskip  $\centerdot \,\,\varphi: S^\bot \to \widetilde S$ : \,\,a projection onto an anticanonical rational surface. \medskip
  
\item We construct in section \textbf{3} the projective surfaces $S$ and $S^\bot$, equipped with natural involutions $\tau$ and $\tau^\bot$, as well as $\widetilde S$, the quotient of $S^\bot$ by $\tau^\bot$.
We then prove that any \textit{hyperelliptic $d$-osculating cover} $\pi: (\Gamma, p) \to (X,q)$ factors through $S^\bot$, and projects onto a rational irreducible curve in $\widetilde S$ (\textbf{3.7. \& 3.8.}). An analogous characterization is in order, for $\pi$ to solve the \textbf{NL Schr\"odinger \& 1D Toda} or the \textbf{sine-Gordon case} (\textbf{3.9.}).\medskip
 
\item In section \textbf{4} we fix a complex elliptic curve $(X,q)=(\mathbb{C}/L, 0)$, and give the original motivation for studying \textit{hyperelliptic $d$-osculating covers} of $X$. We start recalling the definition of the \textit{Baker-Akhiezer} function $\psi_D$, associated to the data ($\Gamma, p, \lambda, D)$, where $\Gamma$ is a smooth complex projective curve of positive genus $g$, $\lambda$ a local parameter at $p \in \Gamma$ and $D$ a non-special effective divisor of $\Gamma$. In case $(\Gamma,p)$ is a hyperelliptic curve marked at a Weierstrass point, we give the Its-Matveev (\textbf{I-M}) exact formula for the \textit{KdV} solution associated to $\psi_D$, as a function of infinitely many variables $\{t_{2j\,\textrm{-}1}, \,j\in \mathbb{N}^*\}$. We end up section \textbf{4} proving that any  \textit{hyperelliptic $d$-osculating cover} of $ \mathbb{C}/L$, gives rise to \textit{KdV} solutions $L$-periodic in $t_{2d\textrm{-}1}$.\medskip
 
\item In section \textbf{5} we take up again the algebraic surface set up developped in section \textbf{3},  recalling that any \emph{hyperelliptic $d$-osculating cover} $\pi:(\Gamma,p) \to (X,q)$ factors through an equivariant morphism $\iota^\bot:\Gamma \to \iota^\bot(\Gamma)\subset S^\bot$, before projecting onto the rational irreducible curve $\widetilde \Gamma :=\varphi\big(\iota^\bot(\Gamma)\big)\subset \widetilde S$. The ramification index of $\,\pi\,$ at $\,p\,$ 
  and the degree of $\iota^\bot: \Gamma \to \iota^\bot(\Gamma)\subset S^\bot$, say $\rho$ and $m$, are natural numerical invariants attached to $\pi$. We also define its \textit{type}, $\gamma = (\gamma_i) \in \mathbb{N}^4$, 
   by intersecting $\iota^\bot_*(\Gamma)$ with four suitably chosen exceptional divisors (\textbf{5.2.}). We assume henceforth that $m=1$ and calculate the linear equivalence class of $\Gamma^\perp \subset  S^\perp$. Basic congruences and inequalities for the latter invariants follow (\textbf{5.4.} \& \textbf{5.5.}). For example, the genus of $\Gamma$ satisifies $(2g+1)^2\leq (2d\,$-$\,1)\big(8n+2d\,$-$\,1\big)$. Any \textit{hyperelliptic cover} solving the other three cases also factors through $S^\perp$ and projects onto a rational irreducible curve in $\widetilde S$. Similar congruences and inequalities for their invariants follow as well (\textbf{5.6.}, \textbf{5.7.} \& \textbf{5.8.})\medskip

\item At last,  in section \textbf{6} we focus on $MH_X(n,d,1,1,\gamma) $, the set of of degree-$n$ \emph{hyperelliptic $d$-osculating covers}, of type $\gamma$, not ramified at the marked point and birational to their natural images in $S^\bot$ (i.e.: such that $\rho=m=1$). For any given $ (n,d) \in \mathbb{N}^*\times\mathbb{N}^*$, we find explicit types $\gamma \in \mathbb{N}^4$ satisfying $\gamma^{(2)}=(2d\,$-$\,1)(2n\,$-$\,2)\,$+$\,3$, for which we give an effective construction (leading ultimately to explicit equations) of the corresponding covers. We thus obtain $(d\,$-$\,1)$-dimensional families of arbitrarily high genus marked curves, solving the \textbf{$d$-th KdV case}. A completely analogous constructive approach can be worked out for the other three cases. \\
      
  \end{enumerate}

\noindent  \section{ Jacobians of curves and hyperelliptic $d$-osculating covers}

 \noindent\textbf{2.1.}
	 Let  $\mathbb
P^1$ denote the projective line over $\mathbb C$ and $(X,q) $ the elliptic curve $(\mathbb{C}/L,0)$, where $L$ is a fixed lattice of $\mathbb{C}$. By a curve we will
mean hereafter a complete integral curve over $\mathbb C$, say $\Gamma$, of positive
 arithmetic genus $ g>0$. If $\Gamma$ is smooth,  its Jacobian variety is a complete connected commutative algebraic group of dimension $g$. For a singular irreducible curve of arithmetic genus $g$ instead, the analogous picture decouples into canonically related pieces, as briefly explained hereafter.\\
We have, on the one hand, the moduli space of degree-$0$ invertible sheaves over $\Gamma$, still denoted by $Jac\,\Gamma$ and called the \textit{generalized Jacobian} of $Jac\,\Gamma$. It is a connected commutative algebraic group, canonically identified to  $H^1(\Gamma,O_\Gamma ^{*})$, with tangent space at its origin equal to $H^1(\Gamma,O_\Gamma)$. In particular, it is $g$-dimensional, although not a complete variety any more.
	 
	The latter is related to the Jacobian variety of the smooth model of $\Gamma$. More generally, let $j:\hat{\Gamma} \rightarrow \Gamma$ be any partial desingularization and consider the natural injection $O_\Gamma \rightarrow j_{*}(O^{*}_{\hat{\Gamma}})$, with quotient $N_j $ , a finite support sheaf of abelian groups. From the corresponding exact cohomology sequence we can then extract\\
	 
	 $\qquad \qquad0\rightarrow H^0(\Gamma,N_j) \rightarrow H^1(\Gamma,O^{*}_\Gamma ) \stackrel{j^*}{\rightarrow} H^1(\hat{\Gamma},O^{*}_{\hat{\Gamma}} ) \rightarrow 0$\\
	 
          $\quad \quad$	 or\\
	  
$\qquad \qquad 0\rightarrow H^0(\Gamma,N_j) \rightarrow Jac\,\Gamma \stackrel{j^*}{\rightarrow}Jac\,\hat{\Gamma} \rightarrow 0$.\\
	 
Hence, the homomorphism	 $j^* :Jac\,\Gamma  \rightarrow Jac\,\hat{\Gamma}, \emph{L} \mapsto j^*(\emph{L})$, is surjective, with kernel the affine algebraic group $H^0(\Gamma,N_j)$. 

On the other hand, we have the moduli space $W(\Gamma)$, of torsionless, zero Euler characteristic, coherent sheaves over $\Gamma$, also called \emph{compactified Jacobian} of $\Gamma$, on which $Jac\,\Gamma$ acts by tensor product. 
Taking direct images by any partial desingularization $j:\hat{\Gamma} \rightarrow \Gamma$, defines an equivariant embedding $j_* : W(\hat{\Gamma}) \rightarrow W(\Gamma)$, such that $\forall \hat{F} \in W(\hat{\Gamma}), \forall L \in Jac\,\Gamma$, we have the projection formula  $j_{*}(\hat{F}\otimes j^{*}(L)) = j_{*}(\hat{F})\otimes L$.
Hence, a $Jac\,\Gamma$-invariant stratification of $W(\Gamma)$, encoding the web of different partial desingularizations between $\Gamma$ and its smooth model. Let me stress that, up to choosing the marked points, any singular irreducible hyperelliptic curves gives rise to KdV, 1D Toda and  NL Schr\"odinger solutions, parameterized by the compactified Jacobian $W(\Gamma)$ (cf. \cite{S-W}6.).

 For any curve $\Gamma$, let $\Gamma ^0 $
and $Jac\,\Gamma$ denote, respectively, the open subset of
smooth points of $\Gamma$ and its \textit{generalized Jacobian}. Recall that
for any smooth point $p\in\Gamma ^0$, the Abel morphism, $
A_p:\Gamma ^0 \to Jac\,\Gamma$, $ p'\mapsto O_\Gamma
(p'$-$p)$, is an embedding and $A_p(\Gamma ^0)$ generates the whole
jacobian. For
any marked curve $(\Gamma , p)$ as above, and any positive integer
$j$, let us consider the exact sequence of $O_\Gamma $-modules $0\to
O _\Gamma \to O_\Gamma (jp)\to O_{jp} (jp)\to 0$, as well as the
corresponding long exact cohomology sequence :\\

$0\to H^0(\Gamma,O_\Gamma )\to H^0\big(\Gamma,O_\Gamma (jp)\big)\to H^0\big(\Gamma,O_{jp}(jp)\big)\stackrel{\delta }\to H^1(\Gamma,O_\Gamma )\to \ldots ,$\\

where $\delta : H^0\big(\Gamma,O_{jp} (jp)\big) \to H^1(\Gamma,O_\Gamma)$ is the cobord morphism and $H^1(\Gamma,O_\Gamma )$ is
canonically identified with the tangent space to $Jac\,\Gamma $ at $ 0$. 

According to the Weierstrass gap Theorem, for any
$d=1,\ldots ,g: = genus(\Gamma) $, there
exists $0<j < 2g$ such that $\delta \Big(H^0\big(\Gamma,O_{jp}(jp)\big)\Big)$ is a
$d$-dimensional subpace, denoted hereafter by $V_{d,p}$. 

For a
generic point $p$ of $\Gamma $ we have $V_{d,p}= \delta
\Big(H^0\big(\Gamma,O_{dp}(dp)\big)\Big)(i.e.:j=d)$.

 In any case, the above filtration
$\{0\} \subsetneq  V_{1,p}\ldots \subsetneq V_{g,p} = H^1(\Gamma,O_
\Gamma)$ is the, so-called, \textit{flag of hyperosculating spaces} to $A_p(\Gamma )$ at $
0 $. For example, $V_{1,p}$ is equal to $\delta \Big(H^0\big(\Gamma,O_p(p)\big)\Big)$, the tangent to $A_p(\Gamma )$ at $0$.\\

 \textbf{Proposition 2.2. }(\cite{T.2}$1.6.$)\\
\textit{Let $(\Gamma, p, \lambda) $ be a hyperelliptic curve, equipped with a local parameter $\lambda$ at a smooth 
Weierstrass point $p \in \Gamma^0$, and consider, for any odd integer $j=2d\,$-$1\geq1$, the exact sequence of $O_\Gamma$-modules}: 
\begin{displaymath}
 0\to O_\Gamma \to O_\Gamma (jp)\to O_{jp}(jp) \to 0  \quad , 
\end{displaymath}
\textit{as well as its long exact cohomology sequence}
 
\begin{displaymath}
 0\to H^0(\Gamma,O_\Gamma )\to H^0\big(\Gamma,O_\Gamma (jp)\big)\to H^0\big(\Gamma,O_{jp}(jp)\big) \stackrel {\delta }{\to }  H^1(\Gamma,O_\Gamma)\to \dots,
\end{displaymath}

$\delta$  \textit{   being the cobord morphism.}\\

\noindent \textit{For any, $m\geq 1$, we also let $ \, [\lambda^{\textrm{-}\,m}]$ denote the class of $\lambda^{\textrm{-}\,m}$ in $H^0\big(\Gamma,O_{mp}(mp)\big)$.
Then $V_{d,p}$ is generated by $\Big\{ \delta \big([\lambda^{2l\textrm{-}1}]\big),\, l = 1,..,d \Big\}$. In other words, the $d$-th osculating subspace to $A_p(\Gamma )$ at $ 0$ is equal to $\delta\Big(H^0\big(\Gamma,O_{jp}(jp)\big)\Big)$, for} $j=2d\,$-$\,1$.\\

\textbf{Definition 2.3.} \\
 \textit{A finite marked morphism $\pi:(\Gamma,p)\to(X,q)$, such that $\Gamma$ is a hyperelliptic curve and $p \in \Gamma$ a smooth Weierstrass point, will be called a hyperelliptic cover. Let $[\,$-$\,1]:(X,q) \to (X,q)$ denote the canonical symmetry, fixing the origin $q \in X$, as well as the three other half-periods $\{\omega_j ,j=1,2,3\}$, and $\tau_\Gamma :(\Gamma, p) \to (\Gamma,p)$ the hyperelliptic involution. Let us recall that the quotient curve $\Gamma/\tau_\Gamma$ is isomorphic to $\mathbb{P}^1$ and }$[\,$-$\,1] \circ \pi=\pi \circ \tau_\Gamma$.\\

 \textbf{Definition 2.5.}\\
 \textit{Let $\pi:(\Gamma,p)\to(X,q)$ be a finite marked morphism and let
$\iota_\pi:X\to \ Jac\,\Gamma $ denote the group
homomorphism $q'\mapsto A_p\big(\pi^*(q'\textrm{-}\,\,q)\big)$. We will say that $\pi$ has osculating order $d$, or equivalently, that it is a
$d$-osculating cover, if $T_oX\subset H^1(\Gamma,O_\Gamma)$, the tangent to $\iota_\pi(X)$ at $ 0 $ is
contained in $V_{d,p}$\textbackslash $V_{d\textrm{-}1,p}$. If $\pi$ also happens to be a hyperelliptic cover, we will simply say that it is a hyperelliptic $d$-osculating cover.} \vspace*{2mm} 
 
\noindent

 The \textit{osculating order} of $\pi$ is a geometrical invariant, bounded by the arithmetic genus of $\Gamma$, which we may want to know. The following hyperelliptic $d$-osculating criterion, analog to Krichever's tangential one
(cf. \cite{K.2} p.289), will be instrumental for its calculation, as well as for further developpment in section \textbf{5.}.\\

\textbf{ Theorem 2.6.}

 \noindent \emph{Let $\pi :(\Gamma ,p)\rightarrow (X,q)$ be an arbitrary hyperelliptic cover of arithmetic genus $g$. Then its osculating order  $d\in \{1,\dots,g\}$ is characterized by the existence of a projection $%
\kappa :\Gamma \rightarrow \mathbb{P}^1$ such that:}

\medskip\ (1) \emph{the poles of $\kappa $ lie along $\pi ^{\textrm{-}1}(q)$};

\medskip\ (2) $\kappa $+$\,\pi ^{*}(z^{\,\textrm{-}1})$ \emph{has a pole of order $2d\,$-$1$  at $p$, and no
other pole along $\pi ^{\textrm{-}1}(q)$}.\\

\noindent \emph{Furthermore, if }$\tau _\Gamma :\Gamma \rightarrow \Gamma $ \emph{denotes the
hyperelliptic involution of $\Gamma $, there exists a unique projection} $%
\kappa :\Gamma \rightarrow \mathbb{P}^1$ \emph{satisfying properties} (1) \& (2)
\emph{above, as well as}:

\medskip\ (3) $\tau _\Gamma ^{*}(\kappa )=\,$-$\,\kappa $.\\

\textbf{Proof}. According to \textbf{2.2.}, $ \forall k\in \{1,..,g\}$ the $k$-th osculating subspace $V_{k,p}$ is generated by 
$\Big\{\delta \big([\lambda^{\textrm{-}(2l\textrm{-}1)}]\big), l = 1,..,k\Big\}$. On the other hand, the 
 tangent to $\iota_\pi(X)\subset Jac\,\Gamma$ at $0$ is equal to $\pi^*\big(H^1(X,O_X)\big)$ and generated by $\delta\big([\pi^*(z^{\textrm{-}1})]\big)$. In other words, the \textit{osculating order} $d$ is the smallest positive integer such that $\delta\big([\pi^*(z^{\textrm{-}1})]\big)$ is a linear combination $\sum^{d}_{l=1}a_l\delta\big([\lambda^{\textrm{-}(2l\textrm{-}1)}]\big)$, with $a_d\neq 0$. Or equivalently, thanks to the Mittag-Leffler Theorem, if and only if there exists a projection $\kappa: \Gamma \to \mathbb{P}^1$, with polar parts equal to 
 $\pi^*(z^{\,\textrm{-}1})\,\textrm{-}\sum^{d}_{l=1}a_l\lambda^{\textrm{-}(2l\textrm{-}1)}$. The latter conditions on $\kappa$ 
 are equivalent to \textbf{2.6.}(1) \& (2). Moreover, up to replacing $\kappa$ by 
 $\frac{1}{2}\big(\kappa \,$-$\, \tau^*_\Gamma (\kappa)\big)$, we can assume $\kappa$ is 
 $\tau_\Gamma$-anti-invariant. Now, the difference of two such functions should be 
 $\tau_\Gamma$-anti-invariant, while having a unique pole at $p$, of order strictly smaller than $2d\,$-$1\leq  2g\,$-$1$.
 But the latter functions are all $\tau_\Gamma$-invariant, implying that the projection $\kappa$ (
 satisfying conditions \textbf{2.6.}(1), (2) \& (3) ), is unique. $\blacksquare$\\

\textbf{ Definition 2.7.}

\emph{The pair of marked projections $(\pi, \kappa)$, satisfying} \textbf{2.6.}(1), (2) \& (3), \emph{will be called  a hyperelliptic d-osculating pair, and $\kappa$ the hyperelliptic $d$-osculating
function associated to $\pi $. In the latter case, $\pi$ gives rise to solutions of the KdV hierarchy, $L$ periodic in the $d$-th KdV flow, as will be proved in section} \textbf{4.}\\

The following Proposition calculates the tangent at any point of the curve $A_{p}(\Gamma)\subset Jac\,\Gamma$, and leads to a useful characterization of the hyperelliptic covers solving the other cases. Its proof follows along the same lines as \textbf{2.2.}'s proof.\\

\textbf{Proposition 2.8.} 

\emph{Let $(\Gamma,r,\lambda )$ be a hyperelliptic curve equipped with a local parameter at an arbitrary smooth point $r\in \Gamma$. Then $V_{\Gamma,r}^1\subset H^1(\Gamma,O_\Gamma)$, the tangent line to $A_{p}(\Gamma)$ at $A_{p}(r)$, is generated by $\delta\big([\lambda^{-1}]\big)$}.\\

\textbf{Corollary 2.9.}

\emph{Let $\pi :(\Gamma ,p)\rightarrow (X,q)$ be an arbitrary hyperelliptic cover, $p^+\in \Gamma$ a non-Weierstrass point, $p^-:=\tau_\Gamma(p^+)$, and let $T_oX\subset H^1(\Gamma,O_\Gamma)$ denote the tangent line defined by the elliptic curve $\iota_\pi(X)\subset Jac\,\Gamma$. Then, the data $(\pi,\,p^+,p^-)$ solves the \textbf{NL Schr\"odinger \& 1D Toda case} (i.e.: $T_oX=V_{\Gamma,p^+}^1=V_{\Gamma,p^-}^1$), if and only if there exists a projection $%
\kappa :\Gamma \rightarrow \mathbb{P}^1$ such that:}

\medskip\ (1) \emph{the poles of $\kappa $ lie in $\pi ^{\textrm{-}1}(q) \cup \{p^+,p^-\}$}.\\

\medskip\ (2) $\kappa $+$\,\pi ^{*}(z^{\,\textrm{-}1})$ \emph{has simple poles at $\{p^+,p^-\}$, and no
other pole along $\pi ^{\textrm{-}1}(q)$}.\\

\medskip\ (3) $\tau _\Gamma ^{*}(\kappa )=\,$-$\,\kappa $.\\ 

\textbf{Corollary 2.10.}

\emph{Let $\pi :(\Gamma ,p)\rightarrow (X,q)$ be an arbitrary hyperelliptic cover equipped with two Weierstrass points $p_o,p_1$, and let $T_oX\subset H^1(\Gamma,O_\Gamma)$ denote the tangent line defined by the elliptic curve }$\iota_\pi(X)\subset Jac\,\Gamma$\emph{. Then, the data $(\pi,\,p_o,p_1)$ solves the \textbf{sine-Gordon case} (i.e.: $T_o X\subset V_{\Gamma,p_o}^1 + V_{\Gamma,p_1}^1$), if and only if there exists a projection $%
\kappa :\Gamma \rightarrow \mathbb{P}^1$ such that:}

\medskip\ (1) \emph{the poles of $\kappa $ lie in $\pi ^{\textrm{-}1}(q) \cup \{p_o,p_1\}$}.\\

\medskip\ (2) $\kappa $+$\,\pi ^{*}(z^{\,\textrm{-}1})$ \emph{has simple poles at $\{p_1,p_2\}$, and no
other pole along $\pi ^{\textrm{-}1}(q)$}.\\

\medskip\ (3) $\tau _\Gamma ^{*}(\kappa )=\,$-$\,\kappa $.\\ 

\section{The algebraic surface set up}

\textbf{3.1.} We will construct hereafter a ruled surface $\pi_S:S \to X$, as well as a blowing-up $e:S^\perp \to S$, having a natural involution $\tau^\perp : S^\perp \to S^\perp$, such that any \emph{hyperelliptic osculating cover}
  $\pi: (\Gamma,p) \rightarrow (X,q)$ factors through $\pi_{S^\bot}$,
   via an equivariant morphism $\iota^\bot:\Gamma \to  \Gamma^\bot:=\iota^\bot(\Gamma)\subset S^\bot$ (i.e.:
    $\iota^\bot \circ \,\tau_\Gamma = \tau^\bot \circ \,\iota^\bot$). We will also prove that $\widetilde \Gamma := \varphi(\Gamma^\bot)$, its image in the quotient surface $\widetilde S: =S^\perp/ \tau^\perp$, is
  an irreducible rational curve. Generally speaking, our main strategy, fully developped in section \textbf{5.}, will consist in translating numerical invariants of $\pi: (\Gamma,p) \to (X,q)$, in terms of the numerical equivalence class of the corresponding rational irreducible curve $\widetilde \Gamma \subset \widetilde S$ and its geometric properties. \\

     \indent The whole relationship is sketched in the diagram below.

\begin{displaymath}
\xymatrix{
&\Gamma^\bot  \subset S^\perp \ar@/_/[ddr]|{\pi_{S^\bot}} \ar[dr]^{e} \ar[r]^{\varphi} &   \widetilde \Gamma \subset\widetilde S\\
p\in\Gamma \ar[rrd]^\pi \ar[ru]^{\iota^\bot} &&S\ar[d]|{\pi_{S}} &  \\
& &q \in X  
     }
\end{displaymath} 
  \textbf{ Definition 3.3.}
 \begin{enumerate}
 \item \textit{Besides the origin $\omega_o:=q\in X$, there are three other half-periods, say $\{\omega_1,\omega_2,\omega_3\}\subset X$, fixed by the canonical symmetry $ [\,$-$\,1]:(X,q) \to (X,q)$}.\\
 
\item \emph{Consider the open affine subsets $U_o:=X\setminus \{q\}$ and $U_1:=X\setminus \{\omega_1\}$ and fix an odd meromorphic function $\zeta:X \to \mathbb{P}^1$, with divisor of poles  equal to $(\zeta)=q+\omega_1\,$-$\,\omega_2\,$-$\,\omega_3$. Let $\pi_S: S \to X$ denote the ruled surface obtained by identifying $\,\mathbb{P}^1\times U_o\,$ with  $\,\mathbb{P}^1\times U_1\,$ over $X \setminus \{q,\omega_1\}$ as follows}:\\

\indent $\forall q' \neq q,\omega_1,\quad (T_o\,,q') \in  \mathbb{P}^1\times U_o\,\,\,$ \emph{is glued with} $\,\,\,(T_1+ \frac{1}{\zeta(q')}\,,q') \in  \mathbb{P}^1\times U_1$.\\

\emph{In other words, we glue the fibers of $\,\mathbb{P}^1\times U_0\,$ and $\,\mathbb{P}^1\times U_1\,$, over any $q' \neq q,\omega_1$, by means of a translation. In particular the constant sections $q' \in U_i \mapsto (\infty,q')\in \mathbb{P}^1\times U_i$, for $i\in\{0,1\}$, get glued together, defining a particular one of zero self-intersection, denoted by $C_o \subset S$}. \\

\item \textit{The meromorphic differentials $dT_o$ and $dT_1$ get also glued together, implying that $K_S$, the canonical divisor of $S$ is represented by -$2C_o$. Any section of $\pi_S: S \to X$, other than $C_o$, is given by two non-constant morphisms $f_i:U_i \to \mathbb{P}^1$ ($i=1,2$), such that $ \quad f_o=f_1\,$-$\,\frac{1}{\zeta}\quad$ outside $\{q,\omega_1\}$. 
A straightforward calculation shows that any such a section intersects $C_o$, while having self-intersection number greater or equal to $2$. It follows from the general Theory of Ruled Surfaces (cf. \cite{H}V.2) that $C_o$ must be the unique section with zero self-intersection}.\\

\item \textit{The only irreducible curve linearly equivalent to a multiple of $C_o$ is $C_o$ itself \big(cf. }\cite{T-V.1}3.2.(1)\big).\\

\item \emph{The involutions $\,\mathbb{P}^1\times U_i \to\,\mathbb{P}^1\times U_i,\quad (T_i,q') \mapsto \big(\,$-$\,T_i,[\,$-$\,1](q')\big)\,\,\,(i=0,1)$, get glued under the above identification and define the involution $\tau: S \to S$, such that $\pi_S \circ \tau=[\,$-$\,1] \circ \pi_S$, already mentioned in} \textbf{3.1.}. \emph{In particular, $\tau$ has two fixed points over each half-period $\omega_i$, one in $C_o$, denoted by $s_i$, and the other one denoted by $r_i$ ($i=0,..,3$)}.\\

\item L\textit{et $e :  S^\bot \to S $ denote hereafter the blow-up of $S$ at
$\{s_i, r_i,i=0,..,3\}$, the eight fixed points of $ \tau$, and $\tau
^\bot : \ S^\bot \to S^\bot$ its lift to an involution fixing the
corresponding exceptional divisors $\big\{s_i^\bot:=e^{-1}(s_i), r_i^\bot:=e^{-1}(r_i), i=0,..,3\big\}$.
Taking the quotient of $S^\bot$ with respect to $\tau
^\bot$, we obtain a degree-$2$ projection $\varphi: S^\bot \to \widetilde S $ onto a smooth rational surface $\widetilde S $, ramified along the exceptional curves $\{s_i^\bot, r_i^\bot, i=0,..,3\}$.
 Let $C_o^\bot$ and $\widetilde C_o$ denote, respectively, the strict transform in $S^\bot$ of $C_o\subset S$ (respectively: the corresponding projections in $\widetilde S$). 
 For any $i=0,..,3$, let also $\widetilde s_i$ and $\widetilde r_i$ denote the projections in $\widetilde S$ of $s_i^\bot$ and $r_i^\bot$, respectively. The canonical divisor of $\widetilde S$, say  $\widetilde K$, satisfies $\varphi^*(\widetilde K)= e^*(\,$-$\,2C_o)$ and is linearly equivalent to -$\,2 \widetilde C_o\,$-$\,\sum_{i=0}^3 \widetilde s_i$. }
\\

\end{enumerate}

  The Lemma and Propositions hereafter, proved in \cite{T.2}$ 2.3., \,2.4. \,\& \,2.5.$, will be 
instrumental in constructing the equivariant factorization $\iota^\bot : \Gamma \to S^\bot$ (\textbf{3.1.}).\\ 

  \textbf{Lemma 3.6.} 

\noindent \textit{There exists a unique, $\tau$-anti-invariant, rational morphism $\kappa_s : S \to \mathbb{P}^1$,  
with poles over $C_o$}+\,\textit{$\pi_S^{\textrm{-}1}(q)$, such that over a suitable neighborhood 
$U$ of $q\in X$, the divisor of poles of $\kappa_s $}+\textit{\,$ \pi_S^*(z^{\textrm{-}1})$ is reduced and equal to
 $C_o\cap \pi_S^{\textrm{-}1}(U)$}. \\

  \textbf{Proposition 3.7.} 
\noindent \textit{For any hyperelliptic cover $\pi: (\Gamma, p) \rightarrow (X, q)$, the 
following conditions are equivalent:}
\begin{enumerate}
\item \textit{there is a projection $\kappa :\Gamma \rightarrow \mathbb{P}^1$, satisfying properties} \textbf{2.6.}(1), (2) \& (3) ;\\

\item \textit{there is a morphism $\iota: \Gamma \to S$ such that $\,\pi = \pi_S \circ \, \iota\,$, $\,\iota\circ \tau_\Gamma=\tau \circ \iota$ and 
$\iota^*(C_o)=(2d\,$-$1)p$.} 
\end {enumerate}
\textit{In the latter case, $\pi$ is a hyperelliptic $d$-osculating morphism (}\textbf{2.5.}) and solves the \textbf{$d$-th KdV case}.\\

 \textbf{Proposition 3.8.}\\
\noindent \textit{For any hyperelliptic $d$-osculating pair $(\pi,\kappa)$, the above morphism
$\iota: \Gamma \to S$ lifts to a unique equivariant morphism $\iota^\bot: \Gamma \to S^\bot$ 
(i.e.: $\tau^\bot \circ \iota^\bot = \iota^\bot \circ \tau_{\Gamma})$. 
In particular, $(\pi,\kappa)$ is the pullback of $(\pi_{S^\bot},\kappa_{s^\bot})=
(\pi_{S}\circ e,\kappa_{s}\circ e)$, and $\Gamma$ lifts to a $\tau^\bot$-invariant curve, $\Gamma^\bot:=\iota^\bot(\Gamma)\subset S^\bot$, which projects onto the rational irreducible curve $\widetilde \Gamma:=\varphi\big(\Gamma^\bot\big) \subset \widetilde S$.}
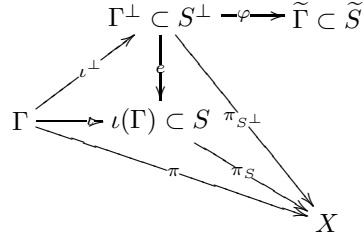
\begin{displaymath}
\xymatrix{
&  \Gamma^\bot \subset  S^\bot \ar[r]|\varphi \ar[d]|e \ar[ddr]|{\pi_{S^\bot}} & \widetilde \Gamma \subset   \widetilde S  \\
\Gamma \ar[r]|\iota \ar[ru]|{\iota^\bot} \ar[rrd]|\pi & \iota (\Gamma) \subset S \ar[rd]|{\pi_S} & \\
& & X
     }
\end{displaymath} 
 
 \textbf{Proof}.  The blow-up $e : S^\bot \to S$, as well as $\iota : \Gamma \to S$, can be pushed down to the corresponding quotients, making up the following diagram:

 
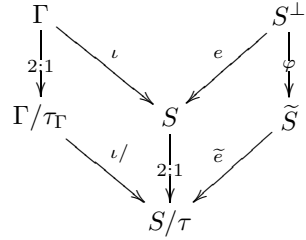
\begin{displaymath}
\xymatrix{
\Gamma  \ar[d]|{2:1} \ar[rd]^\iota &    &    S^\bot \ar[d]|\varphi \ar[ld]_e  \\
\Gamma/\tau_\Gamma  \ar[rd]^{\iota/} & S \ar[d]|{2:1}   &    \widetilde S \ar[ld]_{\widetilde e}  \\
                                                & S/\tau   &   }
\end{displaymath}

Moreover, since $\widetilde e : \widetilde S \to S/\tau$ is a birational morphism and $\Gamma/\tau_\Gamma$ is a smooth curve (in fact isomorphic to $\mathbb{P}^1\,$), we can lift  $ \iota/ :\Gamma/\tau_\Gamma \to S/\tau$  to  $\widetilde S$, obtaining a morphism    $\widetilde \iota:\Gamma   \to \widetilde S$, fitting in the diagram:
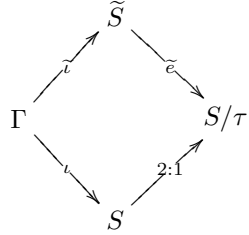
\begin{displaymath}
\xymatrix{
&  \widetilde S \ar[rd]|{\widetilde e} & \\
\Gamma \ar[ru]|{\widetilde \iota} \ar[rd]|\iota & & S/\tau \\
& S \ar[ru]|{2:1} &
  }
\end{displaymath}
 Recall now that $S^\bot$ is the fibre product of  $\,\widetilde e: \widetilde S \to S/\tau$  and  $ S \to S/\tau $ (cf. \cite{T-V.1}4.1.). Hence, $\,\iota\,$  and  $\,\widetilde \iota\,$  lift to a unique equivariant morphism $\iota^\bot : \Gamma \to S^\bot$, fitting in 

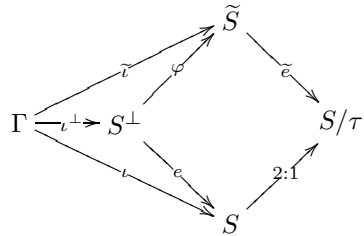
\begin{displaymath}
\xymatrix{
& &  \widetilde S \ar[rd]|{\widetilde e}\\
\Gamma \ar[rru]|{\widetilde \iota} \ar[rrd]|\iota \ar[r]|{\iota^\bot} & S^\bot \ar[ru]|\varphi \ar[rd]|e & & S/\tau\\
& &  S \ar[ru]|{2:1} }
\end{displaymath}

Furthermore, since $\widetilde \iota: \Gamma \rightarrow \widetilde S$ factors through $\Gamma \rightarrow \Gamma /\tau_{\Gamma}\cong \mathbb{P}^1$, its image  $\widetilde \Gamma:=\varphi\big(\iota^\bot(\Gamma)\big)=\widetilde {\iota} (\Gamma)  \subset \widetilde S$ is a rational irreducible curve as claimed.$\blacksquare$\\

Analogously to the \textbf{KdV case}, any data $(\pi,p^+,p^-)$ or $(\pi,p_1,p_2)$, solving the 
\textbf{NL Schr\"odinger \& 1D Toda} or the \textbf{sine-Gordon case}, factors through an equivariant morphism $\iota^\perp: \Gamma \to S^\perp$, and its image $\Gamma^\perp:= \iota^\perp(\Gamma)$ projects onto a rational irreducible curve in $\widetilde S$.\\

 \textbf{Proposition 3.9.} 
\noindent \textit{Let $\pi: (\Gamma, p) \rightarrow (X, q)$ be an arbitrary hyperelliptic cover equipped with two points $p'\neq p'' \in \Gamma$ such that the (divisor) sum $p'+p''$ is $\tau_\Gamma$-invariant. Then, the
following conditions are equivalent:}
\begin{enumerate}
\item \textit{there is a projection $\kappa :\Gamma \rightarrow \mathbb{P}^1$, satisfying properties} \textbf{2.9.}(1), (2) \& (3) \textit{or} \textbf{2.10.}(1), (2) \& (3);\\

\item \textit{there is a morphism $\iota: \Gamma \to S$ such that $\,\pi = \pi_S \circ \, \iota\,$, $\,\iota\circ \tau_\Gamma=\tau \circ \iota$ and 
$\iota^*(C_o)=p' +p''$.} 
\end {enumerate}
\textit{In the latter case, $(\pi,p',p'')$ solves, either the }\textbf{NL Schr\"odinger \& 1D Toda case}, if $\tau_\Gamma(p')=p''$, or the \textbf{sine-Gordon case}, if $p'$ and $p''$ are Weierstrass points.\\

\textbf{Proposition 3.10.}\\
\noindent \textit{For any data $(\pi,p',p'',\kappa)$ as in }\textbf{3.9.}\textit{, the morphism
$\iota: \Gamma \to S$ lifts to a unique equivariant morphism $\iota^\bot: \Gamma \to S^\bot$ 
(i.e.: $\tau^\bot \circ \iota^\bot = \iota^\bot \circ \tau_{\Gamma})$. 
In particular, $(\pi,\kappa)$ is the pullback of $(\pi_{S^\bot},\kappa_{s^\bot})=
(\pi_{S}\circ e,\kappa_{s}\circ e)$, and $\Gamma$ lifts to a $\tau^\bot$-invariant curve, $\Gamma^\bot:=\iota^\bot(\Gamma)\subset S^\bot$, which projects onto the rational irreducible curve $\widetilde \Gamma:=\varphi\big(\Gamma^\bot\big) \subset \widetilde S$.}\\

\section{Complex hyperelliptic curves and elliptic KdV solitons}                                                                             

\noindent\textbf{4.1.} - Let $\Gamma $ be a smooth complex projective curve of
positive genus $g$, equipped with a local coordinate at $p\in \Gamma $, say $\lambda $,  as well as a
non-special degree-$g$ effective divisor $D$ with support disjoint from $p$. Then the so-called \textit{%
Baker-Akhiezer function} associated to the spectral data $(\Gamma
,p,\lambda ,D)$ and denoted by $\psi _D$, is the unique meromorphic function on $\mathbb{C}^\infty\times \left( \Gamma \backslash\{p\}\right)$ such that for any $\vec t = (t_1, t_2,\dots) \in \mathbb{C}^\infty$:

\begin {enumerate}
\item
 the divisor of poles of $\psi_D(\vec t, )$, on $ \Gamma \backslash\{p\}$, is bounded by $D$;

\item in a neighbourhood of $p$, $\psi_D(\vec t,\lambda)$ has an essential
singularity of type:\\

\begin{displaymath}
\psi _D(\vec t,\lambda )=exp\,\big(\sum^{\infty}_{0<i}t_i\lambda^{\textrm{-}i}\big)\big(1\textrm{+} \sum^{\infty}_{0<i} \xi^D_i(\vec t)\lambda^i\,\big) .
\end{displaymath}
\end {enumerate}
For any $i \geq 1$, differentiating $\psi _D$, either with respect to $t_i$, or $i$ times with respect to $x:=t_1$, we obtain a meromorphic
function with  divisor of poles $D + \,ip\,$ and same type of essential singularity
at $p$ as $\psi _D$. We can therefore construct a differential polynomial of degree $i$ in $\frac{\partial}{\partial x}$, with functions of $\vec t$ as coefficients, say $P_i(\frac{\partial}{\partial x})$, such that $\frac{\partial}{\partial t_i}\psi _D \,$-$\, P_i({\frac{\partial}{\partial x}})\psi _D$ has the same properties as $\psi_D$. The uniqueness of the latter \textit{BA} function implies that $\psi _D(\vec t,\lambda )$ satisfies the (so-called \textit{KP}) hierarchy of partial derivatives equations $\frac{\partial}{\partial t_i}\psi_D= P_i(\frac{\partial}{\partial x})\psi_D $, $i \in \mathbb{N}^*$.

\medskip\textbf{4.2.} Let us suppose in the sequel that $(\Gamma ,p)$ is a hyperelliptic curve, marked
at a Weierstrass point, and $\lambda$ an odd local parameter at $p \in \Gamma$. Or in other words, that there exists a degree-2 projection $f:\Gamma \rightarrow
\mathbb{P}^1$, with a double pole at $p$, and $f(\lambda)= \frac{1}{\lambda^2}\,$+$\,O(\lambda^2)$. It is classically known then that the \textit{BA} function $\psi _D(\vec t,\lambda )$, corresponding to any non-special degree-$g$ effective divisor $D$ of $\Gamma $, does not depend, up to an exponential, on the even variables $\{t_{2j}, \, j\in \mathbb{N}^* \}$. 
\medskip For example, choosing $\lambda$ such that $f(\lambda)=\frac{1}{\lambda^2}$, we will have $\psi _D = exp\,\big(\sum_j t_{2j}f^j\big)\psi_{D\big|\{t_{2j}=0\}}$. 

\medskip It then follows that $\psi_{D\big|\{t_{2j}=0\}}$ solves the \textit{KdV} hierarchy and $u:=\, $-$\,2\frac{\partial}{\partial x} \xi^D_1$ the \textit{Korteweg-deVries} equation:
\begin{displaymath}
\label{eq:KdV} u_{t_3}= \frac{1}{4}(6u\cdot u_x + u_{xxx}) \quad(x:=t_1).
\end{displaymath}

 A more concrete formula, (due to A.Its and V.Matveev, cf. \cite{I-M}), is in order:

\begin{displaymath}
\textbf{(I-M)} \,\,\,\, u(t_1,t_3,t_5,\dots)=2\frac{\partial^2}{\partial x^2}\Big( \log \theta _\Gamma
(Z\, \textrm{-} \sum^{\infty}_{0<j}t_{2j\,\textrm{-}1}U_j\,)\Big) +c\text{,} 
\end{displaymath}
where

i) $\theta _\Gamma :\mathbb{C}^g\rightarrow \mathbb{C}$
 denotes the Riemann theta-function of $\Gamma $;

ii) $Z\in \mathbb{C}^g$ projects onto  $A_p(D)$ and $c\in \mathbb{C}$;

iii) $\forall j \geq 1$,  $(2j)!\cdot U_j=A_p^{(2j\,\textrm{-}1)}(\lambda)_{\big|\lambda=0}$, the \textit{$(2j$-$1)$}-th derivative of $A_p(\lambda)$ at $\lambda=0$. 

\medskip\textbf{Remark 4.3.}
\begin{enumerate}
\item The vectors $\{U_k, 1\leq k \leq j\}$ generate $V_{j,p}$, the $j\,$-th \textit{hyperosculating space} to $A_p(\Gamma)$ at $A_p(p)$ (see \textbf{2.1.}).

\item The above construction of \textit{KdV} solutions can be
generalized to any singular marked hyperelliptic curve $(\Gamma, p)$, as recalled in \cite{S-W}. The corresponding solutions are then parameterized by $W(\Gamma)$, the \textit{compactified jacobian} of $\Gamma $. Roughly speaking, any $ L \in W(\Gamma)$, in the complement of the theta divisor, corresponds to a non-special degree-$g$ effective divisor, with support at the smooth points of $\Gamma$. Working in the frame of Sato's Grassmannian (cf. \cite{Sato}, \cite{S-W}6.), one can still define an analogous \textit{BA} function, as well as a \textit{KdV} solution. Hence, the highest the arithmetic genus,
the biggest the family of \textit{KdV} solutions. We are thus naturally
led to allow singular marked hyperelliptic curves.

\item According to the \textbf{(I-M)} formula, the \textit{KdV} solution $u=\, $-$\,2\frac{\partial}{\partial x} \xi^D_1$ is a $t_{2d\textrm{-}1}$\textit{-elliptic KdV soliton} (i.e.: doubly periodic in $t_{2d\textrm{-}1}$), if and only if $U_d$ generates an elliptic curve $X \subset Jac\,\Gamma $.  Or in other words, if $(\Gamma, p) \rightarrow (X, q)$ is a smooth \textit{hyperelliptic $d$-osculating cover}.

\item We will actually prove that any \textit{KdV} solution associated to a \textit{hyperelliptic $d$-osculating cover}, is doubly periodic in $t_{2d\textrm{-}1}$, without assuming the above \textbf{(I-M)} formula, or that $\Gamma$ is a smooth curve (see \textbf{4.5.}). The original idea goes back to \cite{K.2}, p.288-289.

\end{enumerate}

\medskip\textbf{Notations 4.4.}

Choose a lattice $L \subset \mathbb{C}$, equipped with a $\mathbb{Z}$-basis $(2\omega_1, 2\omega_2)$, such that the elliptic curve $(X, q)$ is isomorphic to the complex torus $(\mathbb{C}/L, 0)$, and let $\zeta(z):\mathbb{C} \rightarrow \mathbb{P}^1$, denote the $\zeta$-Weierstrass meromorphic function. Recall (cf. \cite{K.2}, p.283) that $\zeta$ is holomorphic outside $L$ and characterized by the following properties:
\begin{displaymath}
\forall z \in \mathbb{C}\setminus L \,\, \left\{\begin{array}{lll}
 \zeta(z) = z^{\textrm{-}1}\textrm{+}\, O(z) \quad,\textnormal{in a neighborhood of 0}\in \mathbb{C},\\ \\
\medskip\zeta(z\textrm{+}\, 2\omega_j) = \zeta(z)\textrm{+}\, \eta_j, \,\,\,\, j=1,2\quad,
\end{array} \right.
\end{displaymath}

 for some $\eta_1 ,\eta_2 \in \mathbb{C}$, satisfying Legendre's relation: $\eta_1\, 2\omega_2$\,\,-\,\,$\eta_2\, 2\omega_1 = 2\pi\sqrt{\textrm{-}1}$.
 
 \medskip\textbf{Proposition 4.5.}
 
 \textit{Let $\pi: (\Gamma, p) \rightarrow (X, q)$ be a genus-$g$, hyperelliptic $d$-osculating cover, $\kappa$ the unique hyperelliptic $d$-osculating function associated to $\pi$, and choose $\lambda$, an odd local parameter at $p$, such that $\kappa $}+\,\textit{$ \pi^*(z^{\textrm{-}1}) = \lambda^{\textrm{-}(2d\textrm{-}1)}$. Then, for any non-special degree-$g$ effective divisor $D$, with support disjoint from $p$, the KdV solution $u=\, $-$\,2\frac{\partial}{\partial x} \xi^D_1$ associated to $(\Gamma, p, \lambda, D)$ (see }\textbf{4.2.}\textit{), is $L$-periodic in $t_{2d\textrm{-}1}$}. \\
 
 \textbf{Proof}. Denote again by $\psi_D(\vec t,\lambda )$ the \textit{BA} function associated to $D$. Recall (see \textbf{2.4.}) that $\kappa$ has poles only over $\pi^{\textrm{-}1}(q)$, and 
 \begin{displaymath}
 \kappa \textrm{+} \,\pi^*\big(\zeta(z)\big)= \kappa \textrm{+}\,\pi^*\big(z^{\textrm{-}1}\textrm{+} \,O(z)\big)= \lambda^{\textrm{-}(2d\,\textrm{-}1)} \textrm{+} \,O(\lambda) 
 \end{displaymath}
 has a pole of order $2d\,$-$1$ at $p$.
 We then prove, coupling the properties of $\zeta$ and $\kappa$, that for $j=1,2$, the function
  \begin{displaymath}
 \phi_j(p')=exp\,\bigg(2\omega_j\Big(\kappa(p')\textrm{+} \,\zeta\big(\pi(p')\big)\Big)\textrm{-} \,\eta_j\pi(p')\bigg) 
  \end{displaymath}
 is well defined and holomorphic all over $ \Gamma \setminus\{p\}$, thanks to Legendre's relations, and has an essential singularity at $p$ of the following type:
 
 \medskip $\phi_j(p')=exp\,\Big(2\omega_j\lambda^{\textrm{-}(2d\textrm{-}1)} $+$\, O(\pi(p')\big)\Big)=exp\,\big(2\omega_j\lambda^{\textrm{-}(2d\textrm{-}1)}\big)\big(1$+$\,O(\lambda)\big)$.
 
  \medskip The main final argument run as follows. The uniqueness of the \textit{BA} function $\psi_D(\vec t,\lambda )$ implies that 
  \begin{displaymath}
  \psi_D(\vec t \textrm{+} 2\omega_j\vec{e}_{2d\textrm{-}1},\lambda )= \phi_j(\lambda)\cdot\psi_D(\vec t,\lambda )\,,
  \end{displaymath}
  where $\vec{e}_{2d\textrm{-}1}= (0,\dots0,1,0\dots) \in \mathbb{C}^\infty$ is the vector having a $1$ at the $(2d$-$1)$-th place and $0$ everywhere else. At last, comparing their developments around $p$ we obtain the following equality:
\begin{displaymath}
\frac{\partial}{\partial x}\xi^D_1(\vec t \textrm{+} 2\omega_j\vec{e}_{2d\textrm{-}1}) = \frac{\partial}{\partial x}\xi^D_1(\vec t)\,\,,\,\,\,\,j=1,2.
   \end{displaymath}
   In other words, the \textit{KdV} solution $u=\, $-$\,2\frac{\partial}{\partial x}\xi^D_1$ associated to the data $(\Gamma, p, \lambda, D)$, is $L$-periodic in $t_{2d\textrm{-}1}$. $\blacksquare$\\
   
 \section{ The hyperelliptic $d$-osculating covers as divisors of a surface}
 
 \noindent\textbf{5.1.}
  Let us consider again the algebraic surface set up constructed in section \textbf{3}, with the equivariant factorization of any hyperelliptic $d$-osculating cover through $S^\bot$, and its projection onto a rational irreducible curve $\widetilde \Gamma \subset \widetilde S$. The corresponding diagram of morphisms, given hereafter, will also be useful for the \textbf{NL Schr\"odinger \& 1D Toda} and \textbf{sine-Gordon} cases.

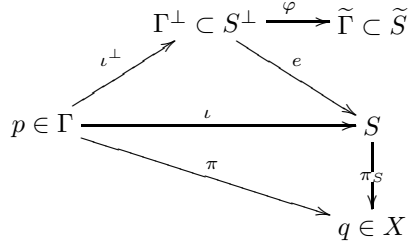
\begin{displaymath}
\xymatrix{
&\Gamma^\bot  \subset S^\perp  \ar[dr]^{e} \ar[r]^{\varphi} &   \widetilde \Gamma \subset\widetilde S\\
p\in\Gamma \ar[rrd]^\pi \ar[ru]^{\iota^\bot} \ar[rr]^\iota &&S\ar[d]|{\pi_{S}} &  \\
& &q \in X  
     }
\end{displaymath} 
 
  \textbf{ Definition 5.2.}\\
\noindent \textit{For any $i=0,..,3$, the intersection number between the divisors $\iota^\bot_*(\Gamma )$ 
and $r_i^\bot$ will be denoted by $\gamma _i$, and the corresponding vector $\gamma = (\gamma_i) \in \mathbb{N}^4$ called the type of $\pi$. Furthermore, $\gamma ^{(1)}$ and
$\gamma^{(2)} $ will denote, respectively, the sums }\\
 
$\quad \quad \quad \quad \quad \gamma ^{(1)} : = \sum^{3}_{i=0}\gamma_i  \qquad and \qquad \gamma ^{(2)} : =  \sum^{3}_{i=0}\gamma_i ^2.$\\

\textbf{Remark 5.3.} \\
\noindent The next step concerns studying the above rational irreducible curves $\widetilde \Gamma \subset \widetilde S$. We will characterize their linear equivalence classes, and dress the basic relations between them and the numerical invariants of the corresponding \textit{hyperelliptic $d$-osculating covers}. These results, already known for $d = 1$ (\cite{T-V.1}) and $d = 2$ (\cite{F}), can be proven within the same framework for any other $d > 2.$\\

 \textbf{Lemma 5.4.} \\
 \noindent \textit{Let $\pi: (\Gamma, p) \rightarrow (X, q)$  be a degree-$n$ hyperelliptic $d$-osculating cover, $\iota^\bot:\Gamma \rightarrow \Gamma^\bot$ its unique equivariant factorization through $S^\bot$ and $\iota := e \circ \iota^\bot$. We let again $\gamma$ denote the type of $\pi$, $\rho$ its ramification index at $p$ and $m$ the degree of $\iota^\bot :\Gamma \to \iota^\bot(\Gamma)$. Then :}

\begin {enumerate}

\item \textit{$ \iota_*(\Gamma )$ is equal to $m.\iota(\Gamma)$ and linearly
equivalent to $nC_o $}+\textit{$\,(2d$-$1)S_o$};

\item \textit{$\iota_*(\Gamma)$ is unibranch, and transverse to the fiber $S_o:=\pi_S^*(q)$ at $s_o= \iota(p)$};

\item \textit{$\rho$ is odd, bounded by $2d\,$-$\,1$ and equal to the multiplicity of $\iota_*(\Gamma)$ at $s_o$};

\item \textit{the degree $m$ divides $n$, $2d\,$-$\,1$ and $\rho$, as well as $\gamma_i, \, \forall i = 0,..,3$};

 \item $\gamma_o $+$ 1 \equiv \gamma_1 \equiv \gamma_2 \equiv \gamma_3 \equiv n(\textnormal{mod}.2)$;
  
  \item \textit{$\iota^\bot_*(\Gamma)$ is linearly equivalent to 
 $e^*\big(nC_o$}+$(2d\,$-$1)S_o\big)\,$-$\,\rho \, s^\bot_o$ -$\,\sum^{3}_{i=0}\gamma_i\, r_i^\bot$.\\
 
\end{enumerate}

\textbf{Proof}. (1) - Checking that $\iota_*(\Gamma)$ is numerically equivalent to $nC_o$+$(2d\,$-$1)S_o$ amounts 
to proving that the intersections numbers $\iota_*(\Gamma) \cdot S_o$ and $\iota_*(\Gamma)\cdot C_o$ are equal 
to $n$ and $2d\,$-$\,1$. The latter numbers are equal, respectively, to the degree of $\pi: \Gamma \rightarrow X$ 
and the degree of $\iota^*(C_o)= (2d\,$-$\,1)p$, hence the result. Finally, since $\iota_*(\Gamma)$ and $C_o$ only intersect at $s_o \in S_o$, we also obtain their linear equivalence.\\
  
 (2) $\,\&\,$ (3) - Let $\kappa: \Gamma \to \mathbb{P}^1$ be the \textit{hyperelliptic $d$-osculating function} associated to $\pi$, uniquely 
characterized by properties \textbf{2.6.}(1), (2) \& (3), and $U \subset X$ a symmetric neighborhood of $q:=\pi(p)$. 
Recall that $\kappa + \pi^*(z^{\textrm{-}1})$ is $\tau_\Gamma$-anti-invariant and well defined over $\pi^{\textrm{-}1}(U)$, and has a (unique) pole of order $2d\,$-$\,1$ at $p$. Studying its trace with respect to $\pi$ 
 we can deduce that $\rho$ must be odd and bounded by $2d\,$-$\,1$.
 
 On the other hand, let $\big(\iota_*(\Gamma), S_o\big)_{s_o}$
  and $\big(\iota_*(\Gamma), C_o\big)_{s_o}$ denote the intersection multiplicities at $s_o$, between 
  $\iota_*(\Gamma)$ and the curves $S_o$ and $C_o$. They are respectively equal, via the projection 
  formula for $\iota \,$, to $\rho$ and $2d\,$-$\,1$. At last, since $\iota_*(\Gamma)$ is unibranch at $s_o$
  and $\big(\iota_*(\Gamma), S_o\big)_{s_o} = \rho \leq 2d\,$-$\,1 = \big(\iota_*(\Gamma), C_o\big)_{s_o}$, 
  we immediately deduce that $\rho$ is the multiplicity of $\iota_*(\Gamma)$ at $s_o$ (and $S_o$ is transverse to $\iota_*(\Gamma)$ at $s_o$).
  
  (4) - By definition of $m$, we clearly have $\iota_*(\Gamma) = m.\iota(\Gamma)$, while 
  $\{\rho, \gamma_i, i=0,..,3\}$ are the multiplicities of $\iota_*(\Gamma)$ at different points of $S$. 
  Hence, $m$ divides $n$ and $2d$-$1$, as well as all integers $\{\rho, \gamma_i, i = 0,..,3\}$.
 
  (5) - For any $i= 0,..,3$, the strict transform of the fiber $S_i:= \pi_S^{-1}(\omega_i)$, by the blow-up 
  $e: S^\bot \to S$, is a $\tau^\bot$-invariant curve, equal to 
  $S^\bot_i:= e^*(S_i) \,$-$\, s_i^\bot \,$-$ \,r_i^\bot$, but also to $\varphi^*({\widetilde S}_i)$, where 
  ${\widetilde S}_i:= \varphi(S^\bot_i)$. Hence, the intersection number 
   $\iota^\bot_*(\Gamma)\cdot S_i^\bot$ is equal to the even integer \\

   $\iota^\bot_*(\Gamma)\cdot S_i^\bot = \iota^\bot_*(\Gamma)\cdot \varphi^*(\widetilde S_i)= 
   \varphi_*(\iota^\bot_*\big(\Gamma)\big)\cdot \widetilde S_i = 2\widetilde \Gamma\cdot \widetilde S_i,$\\
   
   implying that $n= \iota^\bot_*(\Gamma)\cdot e^*(S_i)$ is congruent mod.$2$ to\\
    
 $\iota^\bot_*(\Gamma)\cdot S_i^\bot +\, \iota^\bot_*(\Gamma) \cdot (s_i^\bot +  \,r_i^\bot) \equiv \iota^\bot_*(\Gamma)\cdot (s_i^\bot +  \,r_i^\bot)($mod$.2).$\\
  
   We also know, by definition, that $\gamma_i:= \iota^\bot_*(\Gamma)\cdot r_i^\bot$, while $\iota^\bot_*(\Gamma)\cdot s_o^\bot= \rho$, the multiplicity of $\iota_*(\Gamma)$ at $s_o$, and  $\iota^\bot_*(\Gamma)\cdot s_i^\bot= 0$ if $i \neq 0$, because $s_i \notin \iota(\Gamma)$. Hence, $n$ is congruent mod.$2$, to $\rho $+$\, \gamma_o \equiv 1  $+$ \,\gamma_o\,($mod$.2)$, as well as to $\gamma_i$, if $i\neq 0$.\\
   
  (6) - The Picard group $Pic(S^\bot)$ is the direct sum of $e^*(Pic(S))$ and the rank-$8$
  lattice generated by the exceptional curves $\{s_i^\bot, r_i^\bot, i= 0,..,3\}$. In particular, knowing that $\iota_*(\Gamma)$ is linearly equivalent to $nC_o $+$\, (2d\,$-$\,1)S_o$, and having already calculated  $\iota_*^\bot(\Gamma)\cdot s_i^\bot$ and $\iota_*^\bot(\Gamma)\cdot r_i^\bot$, for any $i= 0,..,3$, we can finally check that $\iota_*^\bot(\Gamma)$ is linearly equivalent to $e^*\big(nC_o $+$(2d\,$-$\,1)S_o\big)\, $-$\,\rho\,  s_o^\bot \,$-$\, \sum^{3}_{0}\gamma_i \, r_i^\bot.  \quad  \blacksquare$ \\
 
 \noindent We are now ready to deduce the basic inequalities relating the numerical
  invariants, associated so far to any such cover $\pi$ (i.e.: $\big\{n, d, g, \rho, m, \gamma\big\}$). The arithmetic genus of the irreducible curve $\widetilde \Gamma:= \varphi(\Gamma^\bot) \subset \widetilde S$, say $\widetilde g$, can be deduced from \textbf{5.4.}(6) via the projection formula for $\varphi : S^\bot \to \widetilde S$. We start proving the inequality $2g+1 \leq \gamma^{(1)}$, before deducing the main one (\textbf{5.5.}(4)) from $\widetilde g \geq 0$. \\

 \textbf{Theorem 5.5.} \\
\noindent \textit{Consider any hyperelliptic $d$-osculating cover $\pi:(\Gamma ,p)\to (X,q)$, of degree $n$, type $\gamma$, arithmetic genus $g$ and ramification index $\rho$ at $p$, and let $m$ denote the degree of its canonical equivariant factorization $\iota^\bot: \Gamma \rightarrow \iota^\bot(\Gamma)\subset S^\bot$. Then the numerical invariants $\{n, d, g, \rho, m, \gamma\}$ satisfy the following inequalities:}

\begin{enumerate}

\item $2g$+$1\leq \gamma ^{(1)}\quad$;\\
\item $\rho=1$ \textit{ implies } $m=1\quad$;\\

\item $\gamma^{(2)} \leq 2(2d\,$-$\,1)(n\,$-$\,m) $+$\, 4\,m^2\,$-$\,\rho^2\quad$;\\

\item $(2g$+$1)^2 \leq 8(2d\,$-$\,1)(n\,$-$\,m) $+$ 13\,m^2$-$\,4\rho^2 \leq 8(2d\,$-$\,1)n $+$ (2d\,$-$\,1)^2.\quad$. \\

\textit{Hence, if $\pi$ is not ramified at $p$, we must have $m=1$, as well as: }\\

\item $(2g$+$1)^2 \leq 8(2d\,$-$1)(n\,$-$\,1)$+$\,9$.\\

\end{enumerate}

\textbf{Proof}. (1) - For any $i=0,..,3$, the fiber of $\pi_{S^\bot}:= \pi_S \circ e: S^\bot \rightarrow X$ over the half-period $\omega_i$, decomposes as $s_i^\bot $+$ \,r_i^\bot $+$\, S_i^\bot$, where $S_i^\bot$ is a $\tau^\bot$-invariant divisor and  $s_i^\bot$ is disjoint with $\iota^\bot_*(\Gamma)$, if $i \neq 0$, while ${\iota^\bot}^*(s_i^\bot)= \rho \, p$, by \textbf{5.4.}(2). Hence, the divisor $R_i:= {\iota^\bot}^*(r_i^\bot)$ of $\Gamma$ is linearly equivalent to $R_i\equiv \pi^{-1}(\omega_i)\,$-$\,(n\,$-$\,\gamma_i)\,p$ (and also $2R_i \equiv 2\gamma_i \,p\,)$. Recalling at last, that $\sum^{3}_{j=1}\omega_j \equiv 3\, \omega_o$, and taking inverse image by $\pi$, we finally obtain that $\sum^3_{i=0}R_i \equiv \gamma^{(1)}\, p\,$. In other words, there exists a well defined meromorphic function, (i.e.: a morphism), from $\Gamma$ to $\mathbb{P}^1$, with a pole of (odd!) degree $\gamma^{(1)}$ at the Weierstrass point $p$. The latter can only happen (by the Riemann-Roch Theorem) if $2g$+$1 \leq \gamma^{(1)}$, as asserted. \medskip

(2) - According to \textbf{5.4.}(4), $m$ divides $\rho$. Hence, $\rho=1$ implies $m=1$.\medskip

(3) - The curve $\iota^\bot(\Gamma)$ is $\tau^\bot$-invariant and linearly equivalent \big(\textbf{5.4.}(4)-(6)\big) to:\\

\quad \quad $\iota^\bot(\Gamma) \sim \frac{1}{m} \Big(e^*\big(nC_o $+$ \,(2d\,$-$1)S_o\big)\,$-$\,\rho  s_o^\bot $-$ \sum ^{3}_{i=0}\gamma_i\,  r_i^\bot\Big)$.\\

Recall also that $\varphi ^{*}(\widetilde{K})$, the inverse image by $\varphi $ of the canonical divisor of $\widetilde{S}$, is linearly equivalent to $\varphi ^{*}(\widetilde{K})\sim e^{*}($-$\,2C_0)$. 
Applying the projection formula for $\varphi : S^\perp \rightarrow \widetilde S$, to the divisor $\iota^\bot(\Gamma)$, we calculate $g(\widetilde \Gamma)$, the arithmetic genus of $ \widetilde \Gamma :=\varphi\,\big(\iota^\bot(\Gamma)\big) \subset \widetilde S$ :\\

\quad \quad $0 \leq g(\widetilde \Gamma) = \frac{1}{4m^2}\Big((2d\,$-$1)(2n\,$-$2m)$+$ \,4m^2$-$\,\rho^2 \,$-$\,\gamma^{(2)} \Big)$, \\

  implying   \\

\quad \quad $\gamma^{(2)} \leq (2d\,$-$1)(2n\,$-$2m)$+$\, 4m^2$-$\,\rho^2$, 

   as claimed. \\

\medskip (4) \& (5) - We start remarking that, for any $j= 1,2,3$, $(\gamma_o\,$-$\,\gamma_j)$ is a non-zero multiple of $m$. Hence,
$\sum_{\substack{i<j}}(\gamma_i\,$-$\,\gamma_j)^2 \geq 3m^2$, and replacing in \textbf{5.5.}(1) we get:
\begin{displaymath}
(2g\textrm{+}1)^2 \leq (\gamma^{(1)})^2 = 4\gamma^{(2)} \,\textrm{-} \,\sum_{\substack{i<j}}(\gamma_i\,\textrm{-}\,\gamma_j)^2 \leq  4\gamma^{(2)}\,\textrm{-}\,\,3m^2.
\end{displaymath}

Taking into account \textbf{5.5.}(3), we obtain the inequality \textbf{5.5.}(4), as well as \textbf{5.5.}(5), which corresponds to the particular case $\rho=m=1$. $\blacksquare$\\ 
 
 \textbf{Lemma 5.6.}\\
  \textit{Let $\pi: (\Gamma, p) \rightarrow (X, q)$ be an arbitrary degree-$n$ hyperelliptic cover, equipped with two points $p'\neq p'' \in \Gamma$ such that the (divisor) sum $p'+p''$ is $\tau_\Gamma$-invariant. Assume  the data $(\pi,\,p',p'')$ solves the \textbf{NL Schr\"odinger \& 1D Toda} or the \textbf{sine-Gordon case}, i.e.: $T_oX=V_{\Gamma,p'}^1+V_{\Gamma,p''}^1$ (}\textbf{2.9. \& 2.10.}\textit{). We let again $\iota:\Gamma \to S$ denote the corresponding morphism (}\textbf{3.10.}\textit{), $\Gamma^\perp$ the image of its lift $\iota^\perp:\Gamma \to S^\perp$, and $\gamma =(\gamma_i)\in \mathbb{N}^4$ its type, obtained by intersecting $\Gamma^\perp$ with the curves $\{r_i^\perp\}$. Then:}

\begin{enumerate}
\item \textit{$\iota(\Gamma)$ is birational to $\Gamma$ and numerically equivalent to $nC_o+2S_o$};\\
\item  \textit{$\iota(\Gamma)$ intersects $C_o$ at $\{\iota(p'),\iota(p'')\}$, with multiplicity $1$ at each point, if 
 $\pi(p')\neq\pi(p'')$, and with multiplicity $2$ if $\pi(p')=\pi(p'')$};\\
 
\item \textit{if $\pi(p')=\pi(p'')$, then $\gamma_o \equiv \gamma_1 \equiv \gamma_2 \equiv \gamma_3 \equiv n(\textnormal{mod}.2)$, $\pi(p')=\omega_{i_o}$ is a half-period and $\iota^\bot_*(\Gamma)$ is linearly equivalent to 
 $e^*\big(nC_o+2S_o\big)\,$-$\,2 s^\bot_{i_o}$ -$\,\sum^{3}_{i=0}\gamma_i\, r_i^\bot$
};\\ 

\item \textit{if $\pi(p')\neq\pi(p'')\notin\{\omega_i\}$, then $\gamma_o \equiv \gamma_1 \equiv \gamma_2 \equiv \gamma_3 \equiv n(\textnormal{mod}.2)$ and $\iota^\bot_*(\Gamma)$ is linearly equivalent to 
 $e^*\big(nC_o+2S_o\big)\,$-$\,\sum^{3}_{i=0}\gamma_i\, r_i^\bot$};\\

\item \textit{if $\pi(p')\neq\pi(p'')$ are two half-periods of $(X,q)$, say $\{\omega_k,\omega_j\}$, for some $k\neq j$, then $\gamma_k+1 \equiv \gamma_j +1 \equiv \gamma_i \equiv \gamma_l \equiv n(\textnormal{mod}.2)$, where $\{j,k,i,l\}=\{0,1,2,3\}$ and $\iota^\bot_*(\Gamma)$ is linearly equivalent to 
 $e^*\big(nC_o+S_k+S_j\big)\,$-$ \, s^\bot_k$ -$ \, s^\bot_j$-$\,\sum^{3}_{i=0}\gamma_i\, r_i^\bot$};\\
 
 \end{enumerate}
 Analogously to what we proved for the \textbf{$d$-th KdV case} (\textbf{5.5.}), we obtain the following relations between the degree and arithmetic genus of the other cases. \\
 
 \textbf{Theorem 5.7. (NL Schr\"odinger \& 1D Toda case)}\\
\textit{Let $\pi: (\Gamma, p) \rightarrow (X, q)$ be an arbitrary degree-$n$ hyperelliptic cover of arithmetic genus $g$, equipped with two points $p^+ \neq p^- \in \Gamma$ exchanged by the hyperelliptic involution $\tau_\Gamma$. Assume $(\pi,\,p^+,p^-)$ solves the \textbf{NL Schr\"odinger \& 1D Toda case} and let $\gamma \in \mathbb{N}^4$ denote its type (}\textbf{5.6.}\textit{). Then, $\gamma_i \equiv  n(\textnormal{mod}.2)$, for any $i$, and}:\\

\begin{enumerate}
\item $2g+2 \leq \gamma^{(1)}$;\\
\item \textit{$\pi (p^+) \neq \pi (p^-)$ implies $\gamma^{(2)}\leq 4n$, as well as $(g+1)^2 \leq 4n$};\\
\item $\pi(p^+)=\pi(p^-)$ \textit{and $n\equiv 0(\textnormal{mod}.2)$ imply }$\gamma^{(2)}\leq 4n\,$-$\,4$\textit{ and} $(g+1)^2\leq 4n\,$-$\,4$;\\
\item $\pi(p^+)=\pi(p^-)$\textit{ and} $n\equiv 1(\textnormal{mod}.2)$ \textit{imply} $\gamma^{(2)}\leq 4n\,$-$\,8$ \textit{and} $(g+1)^2\leq 4n\,$-$\,8$.\\

\end{enumerate}

 \textbf{Theorem 5.8. (sine-Gordon case)}\\
  \textit{Let $\pi: (\Gamma, p) \rightarrow (X, q)$ be an arbitrary degree-$n$ hyperelliptic cover of arithmetic genus $g$, equipped with two Weierstrass points $p_1,p_2 \in \Gamma$. Assume $(\pi,\,p_1,p_2)$ solves the \textbf{sine-Gordon case} and let $\gamma \in \mathbb{N}^4$ denote its type (}\textbf{5.6.}\textit{). Then:}
\begin{enumerate}
\item $2g\leq \gamma^{(1)}$;\\
\item \textit{$\pi(p_1)\neq\pi(p_2)$ implies $\gamma^{(2)}\leq 4n$, as well as $g^2\leq 4n$};\\
\item \textit{$\pi(p_1)=\pi(p_2)$ and $n\equiv 0(\textnormal{mod}.2)$ imply $\gamma^{(2)}\leq 4n\,$-$\,4$ and $g^2\leq 4n\,$-$\,4$};\\
\item \textit{$\pi(p^+)=\pi(p^-)$ and $n\equiv 1(\textnormal{mod}.2)$ imply $\gamma^{(2)}\leq 4n\,$-$\,8$ and $g^2\leq 4n\,$-$\,8$}.\\

\end{enumerate}

\section{ On hyperelliptic $d$-osculating covers of arbitrary high genus}

 $\textbf{6.1.}$ -  Let $C_o^\bot$ denote the strict transform of $C_o$ in $S^\bot$, ${\widetilde C}_o := \varphi(C_o^\perp)$ its projection in $\widetilde S$ and consider an arbitrary degree-$n$ \emph{hyperelliptic $d$-osculating cover} of \textit{type} $\gamma$, say $\pi:(\Gamma,p) \rightarrow (X,q)$, with ramification index $\rho$ at $p$. We will let $\iota^\perp : \Gamma \rightarrow S^\perp$ denote its unique equivariant factorization through $\pi_{S^\perp} : S^\perp \to X$ (\textbf{5.1.}),  $\Gamma^\perp:=\iota^\bot(\Gamma) $ its image in $S^\perp $ and $\widetilde\Gamma$ the corresponding projection into $\widetilde S$. Recall (\textbf{5.4.} \& \textbf{5.5.}) that the above numerical invariants must satisfy the following restrictions :
 \begin{enumerate}
\item $\rho$ is an odd integer bounded by $2d\,$-$\,1$;
 
\item $\gamma_o $+$ 1 \equiv \gamma_1 \equiv \gamma_2 \equiv \gamma_3 \equiv n($mod$. 2)$.\\

Furthermore,  whenever $m:= deg(\iota^\perp : \Gamma \to \Gamma^\perp)$ is equal to $1$ \big(i.e.: $\Gamma$ is birational to $\Gamma^\bot$\big), $\pi$ can be canonically recovered from $\widetilde \Gamma:=\varphi(\Gamma^\perp)$, and  they all satisfy the following properties:\\

  \item  $\widetilde \Gamma$ is an irreducible rational curve of non-negative arithmetic genus equal to $\widetilde g :=\frac{1}{4}\big((2d\,$-$\,1)(2n\,$-$\,2) $+$\, 4\,\,$-$\,\rho^2$-$\,\gamma^{(2)}\big)\geq0$;

 \item  $\Gamma^\perp$ is linearly equivalent to $e^*\big(nC_o + (2d\,$-$1)S_o\big)\,$-$\,\rho {s_o}^\perp \,$-$\sum^{3}_{i=0}{\gamma_i {r_i}^\perp}$;

 \item $\widetilde \Gamma$ intersects ${\widetilde s}_o:=\varphi({s_o}^\perp)$ at a unique point, where it is unibranch and has multiplicity $\rho$;
 \item $\widetilde \Gamma$ intersects ${\widetilde C}_o$ (at most) at ${\widetilde  p}_o:=\widetilde C_o\cap {\widetilde s}_o$ (i.e.: $\widetilde \Gamma\cap {\widetilde C}_o\subset  \widetilde C_o\cap {\widetilde s}_o$), with multiplicity $\frac{1}{2}(2d\,$-$\,1\,$-$\,\rho)$. In particular, if $\rho= 2d\,$-$\,1$, $\widetilde \Gamma$ and  ${\widetilde C}_o$ are disjoint curves.\\
 
 \end{enumerate}

 \textbf{Definition 6.2.} \\
 \textit{For any $(n, d, \rho, \gamma) \in \mathbb{N}^7$ satisfying the above restrictions, we let $\Lambda (n,d,\rho, \gamma)$ denote the unique element of $Pic(\widetilde S)$ such that $\varphi^*\big(\Lambda(n,d,\rho, \gamma)\big)$ is linearly equivalent to $ e^*\big(nC_o + (2d\,$}-\textit{$1)S_o\big)\,$-$\,\rho {s_o}^\perp \,$-$\sum^{3}_{i=0}{\gamma_i {r_i}^\perp}$, and $MH_X(n,d,\rho,1, \gamma)$ denote the moduli space of degree-$n$} \textit{hyperelliptic $d$-osculating covers} \textit{of type $\gamma$, ramification index $\rho$ at their marked point, and birational to their canonical images in $S^\perp$}.\\
 
 \textbf{Remark 6.3.}

We will restrict to the simpler case where $\rho=1$, $\,\Gamma$ is isomorphic to $\Gamma^\perp $ and $\widetilde \Gamma $ is isomorphic to $\mathbb{P}^1$. In other words, we will focus on degree-$n$ \textit{hyperelliptic $d$-osculating covers} with $\rho=m=1$, and of type $\gamma$ satisfying $\gamma^{(2)}=(2d\,$-$\,1)(2n\,$-$\,2)\,$+$\,3$. We will actually choose $\gamma= (2d\,$-$\,1)\mu+2\varepsilon$, where $\mu$ is an arbitrary $\mu \in \mathbb{N}^4$ satisfying $\mu_o \,$+$\,1\equiv \mu_1\equiv \mu_2\equiv \mu_3($mod$.2)$ and $\varepsilon \in \mathbb{Z}^4$ is equal to $\varepsilon = (d\,$-$1,d\,$-$1,d\,$-$1,0)$. Given such triplet $(n,d,\gamma)$ we give a straightforward construction of $MH_X(n,d,1, 1,\gamma)$ as a $(d\,$-$1)$-dimensional family of curves, embedded in $ S^\bot$ (\textbf{6.9.}). Moreover, it can also be proved that any $\pi\in MH_X(n,d,1, 1,\gamma)$ has a unique birational model in $\,\mathbb{P}^1 \times X$, as a linear combination of $d$ specific polynomials with elliptic coefficients. The same can be done for $2\varepsilon =(d\,$+$1,d\,$-$1,d\,$-$1,d\,$-$1)$ if $d$ is odd, or for $2\varepsilon =(d\,$-$2,d,d,d)$ if $d$ is even; or when permuting and/or changing the signs of their coefficients.\\

\indent We will need the following existence and irreducibility criteria.\\

\textbf{Proposition 6.4.} (\cite{T.2}$3.4$)\\
\textit{Any curve $\Gamma \subset S$ intersecting $C_o $ at a unique smooth point $p\in\Gamma$ is irreducible.}\\

\textbf{Proposition 6.5.}\\
\textit{
 Let $\Gamma ^\bot \subset S^\bot$ be a curve with no irreducible component in $\{r_i^\bot, i=0,..,3\}$, and intersecting $C_o^\bot $ (at most) at a unique smooth point $p^\bot\in\Gamma^\bot$. 
Then $\Gamma ^\bot$ is an irreducible curve.}\\

\textbf{Proof}. The properties satisfied by $\Gamma ^\bot$ assure us that it is the
strict transform of its direct image by
$e:S^\bot \to S$, $\Gamma : = e _*(\Gamma ^\bot)$,  and that the latter does not contain $C_o$. We can also check, that $\Gamma $ is smooth at $p:= e(p^\bot)$ and $\Gamma \cap C_o=\{p\}$. It follows, by $\textbf{6.4.}$,
that ($\Gamma $, as well as its strict transform) $\Gamma^\bot$ is an irreducible curve.$\blacksquare$\\

 \textbf{Proposition 6.6.} (\cite{T-V.1}$6.2.$ )

 \noindent \textit{Any $\alpha  = (\alpha _i)\in \mathbb
N^4$ such that $\,\alpha^{(2)}=2n+1$ is odd gives rise to an exceptional curve of the first kind $\widetilde \Gamma_\alpha \subset \widetilde S$. More precisely, let $k\in \{0,1,2,3\}$ denote the index satisfying $\alpha _k+1 \equiv\alpha _j$\textnormal{(mod.2)}, for any
$j\neq k$, and $S_k:=\pi_S^{-1}(s_k)$, then $\widetilde \Gamma_\alpha\,$ has self-intersection -$\,1$ and $\varphi^*(\widetilde \Gamma_\alpha) \subset S^\bot$ is the unique $\tau^\bot$-invariant irreducible curve linearly equivalent to
$e^*(nC_o+S_k)\,$-$\,s_k^\bot$-$\sum _{i=0}^3\alpha _i r_i^\bot$}. \\

\textbf{Proof}. Let $\Lambda$ denote the unique numerical equivalence class of $\widetilde S$ satisfying $\varphi^*(\Lambda)=e^*(nC_o+S_k)\,$-$\,s_k^\bot$-$\sum _{i=0}^3\alpha _i r_i^\bot$. It has self-intersection $\Lambda \cdot \Lambda=\,$-$\,1$, and $\Lambda\cdot \widetilde K=\,$-$\,1$ as well. It follows that $h^o\big(\widetilde S, O_{\widetilde S}(\Lambda)\big)\geq \chi\big(O_{\widetilde S}(\Lambda)\big)=1$, hence there exists an effective divisor $\widetilde \Gamma \in \big|\Lambda\big|$. Such a divisor is known to be unique and irreducible (\cite{T-V.1}6.2.).
$\blacksquare$\\

\textbf{Corollary 6.7. (\cite{T-V.1})} \\   
\emph{Let $\alpha \in \mathbb{N}^4$ be such that $\alpha _o+1 \equiv\alpha _j$\textnormal{(mod.2)}, $\widetilde\Gamma_\alpha$ the corresponding exceptional curve (see }\textbf{6.6.}\textit{), and $\Gamma_\alpha^\perp :=\varphi^*(\widetilde\Gamma_\alpha)$ its inverse image in $ S^\perp$, marked at its Weierstrass point $p_\alpha:= \Gamma_\alpha^\perp \cap s_o^\perp$. Then, $(\Gamma_\alpha^\perp ,p_\alpha)$ gives rise to $KdV$ solutions, $L$-periodic in $x=t_1$ (the first $KdV$ flow).}\\

\indent The latter corollary will be generalized as follows:
given any $n,d \in \mathbb{N}^*$, we will construct \emph{types} $\gamma =(2d\,$-$\,1)\mu$+$2\varepsilon \in \mathbb{N}^4$, such that $\gamma_o$+$1\equiv \gamma_1\equiv \gamma_2 \equiv \gamma_3$(mod.2) and $\gamma^{(2)}=(2n\,$-$\,2)(2d\,$-$\,1)$+$3$, for which the linear system $|\Lambda(n,d,1,\gamma)|$  (see \textbf{6.2.}) has dimension $d\,$-$\,1$ and a generic element isomorphic to $\mathbb{P}^1$. Hence, they will give rise to $(d\,$-$\,1)$-dimensional families of marked curves solving the \textbf{$d$-th KdV case}.   \\

\textbf{Theorem 6.8.}\\
\emph{Let us fix $\,d\geq2\,$, $\,k\in \{0,1,2,3\}\,$,  and $\mu\in \mathbb{N}^4$ such that $\mu_o+1\equiv
\mu_j \textnormal{(mod.2)}$ (for $j=1,2,3$). Pick any vector $
 \,2\varepsilon=(2\varepsilon_i)\in 2\mathbb{Z}^4\,$, satisfying $(\forall i=0,\ldots ,3)\,$, either}\\

\quad\quad \quad \quad \quad \quad \quad \quad \, \quad  $\vert 2\varepsilon_i\vert = (2d\,$-$\,2)(1\,$-$\,\delta
_{i,k})\,\quad,$

\begin{displaymath}
or \quad\,\, \left\{\begin{array}{lll}
 |2\varepsilon_i| = d\,\,\textrm{-}\,(\textrm{-}1)^{\delta_{i,k}}\quad\quad \textit{if} \quad d \quad\textit{is odd}\quad,\\ \\
\medskip |2\varepsilon_i| = d\,\textrm{-}\,2 \delta_{i,k}\quad\quad \quad \quad\textit{if} \quad d \quad\textit{is even}\quad ,
\end{array} \right.
\end{displaymath}

\emph{as long as $\gamma:= (2d\,$-$1)\mu + 2\varepsilon \in \mathbb{N}^4$, and let $n$ satisfy $\gamma^{(2)} =  (2d\,$-$1)(2n$-$2)+3$. Then $|\,\varphi^*\big(\Lambda(n,d,1,\gamma)\big) |$ contains a
$(d\,$-$\,1)$-dimensional subspace such that its generic element, say $\Gamma^\perp$, satisfies the following properties}:

\begin{enumerate}
\item  $\Gamma^\perp$ \textit{is a} $\tau^\perp$\textit{-invariant smooth irreducible curve of genus }$g
:\ =\frac{1}{2}($-$1$+$\gamma ^{(1)})$;\\

\item \textit{$\Gamma^\bot$ can only intersect $C_o^\bot$ at }$p_o^\bot:=C_o^\bot \cap s_o^\bot$;\\

\item $\varphi(\Gamma^\perp) \subset \widetilde S$ \textit{is isomorphic to } $\mathbb P^1$.\\
  
\end{enumerate}

\textbf{Corollary 6.9.}\\
\noindent \textit{Given $(n,d,\gamma) \in \mathbb{N}^*\times \mathbb{N}^*\times \mathbb{N}^4$ as above, the moduli space $MH_X(n,d,1,1,\gamma)$  (}\textbf{6.2.}\textit{) has dimension $d\,$-$\,1$, and a smooth generic element of genus $g
:\ =\frac{1}{2}($-$1$}$+$ $\gamma ^{(1)})$.\\

  \textbf{Proof} of \textbf{Theorem 6.8.}.

  We will only  work out the case $\gamma:= (2d\,$-$1)\mu + 2\varepsilon$, with $\,\,\varepsilon = (0,\,d\,$-$\,1,\,
d\,$-$\,1, \,d\,$-$\,1)\,\,$. 
\indent For any other choice of $\varepsilon$, the corresponding proof runs along the same lines and will be skipped. In our case, the arithmetic genus $g$ and the degree $n$ satisfy: 
\begin{displaymath}
2g+1= (2d\,\textrm{-}\,1)\mu^{(1)}+6(d\,\textrm{-}\,1)\quad \textrm{and} \quad 2n=(2d \,\textrm{-}1) \mu ^{(2)}
\textrm{+} \,4(d\,\textrm{-}1)(\mu_1 \textrm{+} \mu_2 \textrm{+} \mu_3)\textrm{+}\,6d\,\textrm{-}7.
\end{displaymath}

\indent Consider $\overline \mu:\ =\mu$+$\,(1,1,1,1),\,\, \mu' :\  =
\mu$+$\,(0,2,1,1),\,\, \mu''=  \mu$+$\,(0,0,1,1)$, and let $\overline{Z}^{\,\bot}
,{Z'}^\bot,{Z''} ^\bot \subset S ^\bot $ denote the unique $\tau ^\bot$-invariant curves
linearly equivalent to:\\

\noindent 1) $\overline{Z}^{\,\bot} \,\,\,\sim \,e^*(\overline{m}\,C_o$+$\,S_o)\,$-$\,s_o^\bot\,$-$\sum
_i\overline{\mu}_i r_i ^\bot $, \,\,where $\,\,2\overline{m} \,\,$+$ 1 = \overline {\mu}^{\,(2)}$;

\noindent 2) ${Z'}^\bot \,\,\sim \,e^*(m'C_o$+$\,S_1)\,$-$\,s_1^\bot$-$\sum _i\mu_i' 
r_i^\bot$ , \,where \,\,$2m'\,$+$1={\mu'}^{(2)}$;

\noindent 3) ${Z''}^\bot \sim \,e^*(m''C_o$+$\,S_1)\,$-$\,s_1^\bot$-$\sum
_i\mu_i'' r_i^\bot$ , where $\,2m''$+$1={\mu''}^{(2)}$.\\

\indent Moreover, if $\mu_o\neq 0$ we choose $\underline {\mu}\  = \mu $+$\,
($-$\,1,1,1,1)$ and $2\underline {m}\,$+$1= \underline {\mu}^{(2)}$, and let $\underline {Z}^\bot\subset S ^\bot $ denote the unique  $\tau ^\bot$-invariant curve $ \underline {Z}^\bot \sim e^*(\underline {m}C_o$+$\,S_o)\,$-$\,s_o^\bot \,$-$\sum _i\underline {\mu}_{i} r_i^\bot$.

\indent However, if $\mu_o=0$ we will simply put $\underline {Z}^\bot :\  = \overline{Z}^{\,\bot}
$+$ \,2r_o^\bot$, so that in both cases,  
the divisors $D^\bot_0 : = \overline{Z}^{\,\bot}  $+$  \underline {Z}^\bot $+$ \,2s^\bot _0$ and $
D^\bot_1 : = {Z'}^\bot $+$ {Z''}^\bot $+$ \,2s_1^\bot$  will be linearly equivalent. 
Let us also define,  
\begin{align}
\mu _{(1)} :\  = \mu''=\mu + (0,0,1,1), \notag \\
\mu_{(2)}:\ = \mu + (0,1,0,1),\notag \\
\mu_{(3)}: \ = \mu + (0,1,1,0),\notag
\end{align}
  and let $ Z^\bot_{(k)}  (k=1,2,3)$ be the $\tau^\bot$-invariant curve of
$S^\bot $, linearly equivalent to $e^*(m_{(k)}C_o$+$\,S_k)\,$-$\,s_k^\bot $-$
\sum_i\mu_{(k)i}r_i^\bot$, where $ 2m_{(k)}$+$ 1= \sum_i\mu_{(k)i}^2$.

\indent At last, consider $Z^\bot \sim e^*(mC_o$+$S_o) \,$-$\,
s^\bot_o\,$-$\, \sum_i\mu_ir_i^\bot$, where $2m$+$1 = \sum_i \mu_i^2$ (\textbf{6.2.}).
The $(d$-$1)$-dimensional subspace of $\big|\varphi^*\big(\Lambda(n,d,1,\gamma)\big)\big|$ we are looking for, will be made of all above curves.
 We first remark the following facts :\\

\noindent a) we can check, via the adjunction formula, that any $\tau^\bot$-invariant 
element of $\big|\varphi^*\big(\Lambda(n,d,1,\gamma)\big)\big|$ has arithmetic genus $g:\ =\frac{1}{2}($-$1$+$\,\gamma^{(1)})$, and is the pull-back by 
$ \varphi : S^\bot \to \widetilde S $, of a divisor of zero arithmetic genus of $ \widetilde S$;\\

\noindent b) the following  $d\,$-$\,1$ divisors 

\begin{displaymath}
\Big\{F^\bot_j : = C^\bot_o \textrm{+} \sum _{k=1}^{3}(Z^\bot_{(k)}\textrm{+}2s^\bot_k) \textrm{+}
jD^\bot_o \textrm{+} \,(d\,\textrm{-}2\,\textrm{-}\,j)D^\bot_1, \,\,j=0,...,d\,\textrm{-}2\Big\},
 \end{displaymath}
as well as 
\begin{displaymath}  
 G^ \bot : = Z^\bot \textrm{+} \,(d\,\textrm{-}1)D^\bot_o, 
 \end{displaymath}
are $\tau^\bot$-invariant, belong to $\big|\varphi^*\big(\Lambda(n,d,1,\gamma)\big)\big|$ and have $p^\bot_o :  = C^\bot_o
\cap s_o^\bot$ as their unique common point;\\

\noindent c) the curve $F^\bot_o$ is smooth at $p^\bot_o$, while any other $F^\bot_j$ has multiplicity $1<2j+1< 2d$ at $p^\bot_o$. In particular, they span a $(d\,$-$\,2)$-subspace of $\big|\varphi^*\big(\Lambda(n,d,1,\gamma)\big)\big|$, having a generic element smooth and  
transverse to $s_o^\bot$ at $p^\bot_o$;\\

\noindent d) the curve $G^\bot$ has multiplicity $2d$ at $p_o^\bot$, and no common irreducible component with any $F^\bot_j (\,\forall j=0, \ldots,d\,$-$\,2)$, implying that $\langle \,G^\bot,\,F_j^\bot,\,j=0,..,d\,$-$\,2\,\rangle \subset \big|\varphi^*(\Lambda)\big|$, the $(d\,$-$\,1)$-subspace they span, is fixed component-free;\\

\noindent e) any irreducible curve $\Gamma^\bot \in \langle \,G^\bot,\,F_j^\bot,\,j=0,..,d\,$-$\,2\,\rangle$ projects onto a smooth irreducible curve (isomorphic to $\mathbb{P}^1$). In particular $\Gamma^\bot$ must be smooth outside $\cup_{i=0}^3r_i^\bot$.\\

\noindent f) the curves $G^\bot$ and $F_o^\bot$ have no common point on any $r_i^\bot$ ($i=0,..,3$), implying that $\Gamma^\bot$, the generic element of $\langle \,G^\bot,\,F_j^\bot,\,j=0,..,d\,$-$\,2\,\rangle$, is smooth at any point of $\cup_{i=0}^3r_i^\bot$ and satisfies the announced properties, i.e.:\\

(1) - $\Gamma^\bot$ is $\tau^\bot$-invariant, smooth and satisfies the irreducibility criterion \textbf{6.5.};\\

(2) - $p_o^\bot$ is the unique base point of the linear system and $\Gamma^\bot \cap C_o^\bot= \{p_o^\bot\}$;\\

(3) - its image $\varphi(\Gamma^\bot) \subset \widetilde S$ is irreducible, linearly equivalent to $ \Lambda(n,d,1,\gamma)$ and of arithmetic genus
$\frac{1}{4}\big((2d\,$-$1)(2n$-$2)$+$\,3\,$-$\, \gamma^{(2)}\big)  = 0$; hence, isomorphic to $\mathbb P^1$. $\blacksquare $\\

\textbf{Proof} of \textbf{Corollary 6.9.}.

 The degree-$2$ projection $\varphi: \Gamma^\bot \rightarrow \varphi(\Gamma ^\bot)$ is ramified at $p_o^\bot$ and $\varphi(\Gamma^\bot)$ is isomorphic to $\mathbb{P}^1$. Moreover, $\Gamma^\bot$ is a smooth irreducible curve linearly equivalent to $\big|\varphi^*\big(\Lambda(n,d,1,\gamma)\big)\big|$, of  arithmetic genus 
$g:=\frac{1}{2}(\gamma^{(1)}\,$-$\,1) $. 

In other words, the natural projection $ (\Gamma ^\bot,\,p_o^\bot) \subset (S^\bot ,\,p_o^\bot)\stackrel{\pi_{S^\bot}}{\longrightarrow} (X,q)$ is a smooth degree-$n$ \textit{hyperelliptic $d$-osculating cover} of type $\gamma $,
and genus $g$, such that $(2n$-$2)(2d\,$-1$)$+$\,3=\gamma ^ {(2)} $ and
$2g\,$+$1=\gamma^{(1)}.$ $  \blacksquare $\\
  
 \textbf{Remark 6.10.}\\
\begin{enumerate}
\item The irreducible components of the $d$ generators $\langle \,G^\bot,\,F_j^\bot,\,j=0,..,d\,$-$\,2\,\rangle$ are well known curves, for which one can provide explicit equations in $\mathbb{P}^1\times X$. Hence, any element of $MH_X(n,d,1,1,\gamma)$  is birational to the zero set of a linear combination of $d$ specific degree-$n$ polynomials with coefficients in \textit{K(X)}, the field of meromorphic functions on $X$.\\
\item Effective solutions to the \textbf{NL Schr\"odinger \& 1D Toda} and \textbf{sine-Gordon} cases can also been found through an analogous method. Roughly speaking, we construct infinitely many $1$-dimensional familes of solutions (for both cases), having arbitrary degree $n$, and arbitrary genus $g$. As we shall see, the results differ on whether the pair of marked points have same projection in $X$ or not (and depend on the parity of $n$ as well). The main results are given below (detailed proofs will be given elsewhere).\\

 \end{enumerate}

 \textbf{Proposition 6.11. (NL Schr\"odinger \& 1D Toda restrictions)}\\
 \textit{Let $\pi: (\Gamma, p) \rightarrow (X, q)$ be an arbitrary hyperelliptic cover, equipped with two non-Weierstrass points $p^+ ,p^- \in \Gamma$, such that } $(\pi,p^+,p^-)$\textit{ solves the} \textbf{NL Schr\"odinger \& 1D Toda case.}\textit{Then, the arithmetic genus of }$\Gamma$ \textit{and the degree of }$\pi$, \textit{say $g$ and $n$, satisfy}:\\

\begin{enumerate} 
\item $(g+1)^2\leq 4n\,$-$\,4\quad$ \textit{, if} $\quad \pi(p^+)=\pi(p^-)\quad$ \textit{and} $\quad n\equiv 0$(mod 2);\\
\item $(g+1)^2\leq 4n\,$-$\,8\quad$ \textit{, if} $\quad \pi(p^+)=\pi(p^-)\quad$ \textit{and} $\quad n \equiv 1$(mod 2);\\
\item $(g+1)^2\leq 4n\,\quad \quad$ \textit{, if} $\quad \pi(p^+) \neq \pi(p^-)$.\\
 \end{enumerate}
 
  \textbf{Proposition 6.12. (sine-Gordon restrictions)}\\
 \textit{Let $\pi: (\Gamma, p) \rightarrow (X, q)$ be an arbitrary hyperelliptic cover, equipped with two Weierstrass points $p_o ,p_1 \in \Gamma$, such that } $(\pi,p_o,p_1)$\textit{ solves the} \textbf{sine-Gordon case.}\textit{Then, the arithmetic genus of }$\Gamma$ \textit{and the degree of }$\pi$, \textit{say $g$ and $n$, satisfy}:\\

\begin{enumerate} 
\item $g^2\leq 4n\,$-$\,4\quad \,\,\,\,$ \textit{, if} $\quad \pi(p_o)=\pi(p_1)\quad$ \textit{and} $\quad n\equiv 0$(mod 2);\\
\item $g^2\leq 4n\,$-$\,8\quad \,\,\,\,$ \textit{, if} $\quad \pi(p_o)=\pi(p_1)\quad$ \textit{and} $\quad n \equiv 1$(mod 2);\\
\item $g^2\leq 4n\,$-$2\quad \quad$ \textit{, if} $\quad \pi(p_o) \neq \pi(p_1)$.\\
 \end{enumerate}
 
 Along with the latter restrictions we have the following effective results.\\
 
\textbf{Theorem 6.13. (odd degree NL Schr\"odinger \& 1D Toda case)}\\
 \textit{For any $\alpha \in \mathbb{N}^4$ and $a \in X$ there exists a hyperelliptic cover $\pi: (\Gamma, p) \rightarrow (X, q)$, equipped with two non-Weierstrass points $p^+ ,p^- \in \Gamma$ such that:}\\
 \begin{enumerate}
 \item $\pi(p^+)=a$, $p^+=\tau_\Gamma(p^-)$ and $(\pi,p^+,p^-)$ \textit{solves the} \textbf{NL Schr\"odinger case};\\
 \item $\Gamma$ \textit{has arithmetic genus} $g:= \alpha^{(1)}+1$;\\
\item $deg(\pi)=\alpha^{(2)}+\alpha^{(1)} +1\quad$ \textit{if} $a\notin \{\omega_i\},\quad$ \textit{hence} $\quad \pi(p^+) \neq \pi(p^-) $;\\
\item $deg(\pi)=\alpha^{(2)}+\alpha^{(1)} +3\quad$ \textit{if} $a\in \{\omega_i\},\quad$ \textit{hence} $\quad \pi(p^+) =\pi(p^-) $.\\
\end{enumerate}

\textbf{Theorem 6.14. (even degree NL Schr\"odinger \& 1D Toda case)}\\
 \textit{For any $\alpha \in \mathbb{N}^4\setminus\{0\}$ and $a\in X$ such that, either $\alpha^{(1)}\equiv 0$}(mod 2)\textit{ and $a \notin \{\omega_i\}$, or  $\alpha^{(1)}\equiv 1$}(mod 2)\textit{ and $a \in \{\omega_i\}$, there exists a hyperelliptic cover $\pi: (\Gamma, p) \rightarrow (X, q)$, equipped with two non-Weierstrass points $p^+ ,p^- \in \Gamma$ such that:}\\
 \begin{enumerate}
 \item $\pi(p^+)=a$, $p^+=\tau_\Gamma(p^-)$ \textit{and} $(\pi,p^+,p^-)$ \textit{solves the}\textbf{ NL Schr\"odinger case;}\\
 \item $\Gamma$ \textit{has arithmetic genus} $g:= \alpha^{(1)}\,$-$\,1$;\\
\item $deg(\pi)=\alpha^{(2)}$ \textit{if} $a \notin \{\omega_i\}$\textit{, and} $deg(\pi)=\alpha^{(2)}+1$ \textit{otherwise}.\\

\end{enumerate}

For a better presentation of our sine-Gordon's results, we must also take in account the projections of $(p_o,p_1)$, the pair of Weierstrass points (see \textbf{2.10.}). They either project onto the same point, which can be chosen equal to $\pi(p_o)=\pi(p_1)=\omega_o$, or their projections differ by a non-zero half-period, say $\pi(p_o)=\omega_o$ and $\pi(p_1)=\omega_1$. In all four cases we find $1$-dimensional families of solutions. Additional properties, such as the existence of a fixed point free involution or a real structure can also be found. For example, if $(X,q)$ has a real structure, we can extract from the first three \textbf{sine-Gordon} cases a real $1$-dimensional family having a real structure fixing the Weierstrass points.\\

   \textbf{Theorem 6.15. (even degree sine-Gordon with distinct projections)}\\
   \textit{Pick any $\alpha \in \mathbb{N}^4$ satisfying $\alpha_{2}+\alpha_3\equiv 1$}(mod 2). \textit{Then, there exists a $1$-dimensional family of hyperelliptic covers $\pi: (\Gamma, p) \rightarrow (X, q)$, equipped with a pair of distinct Weierstrass points $\{p_o ,p_1\} \in \Gamma$, such that:}\\
 \begin{enumerate}
 \item $\pi(p_j)=\omega_j$,  \textit{for $j=0,1$ and} $\quad(\pi,p_o,p_1)\,\,$ \textit{solves the} \textbf{sine-Gordon case};\\
 \item $\Gamma$ \textit{has arithmetic genus} $\quad g:= \alpha^{(1)}+1\quad$ \textit{and} $\quad deg(\pi)=\alpha^{(2)}+\alpha_o+\alpha_1 +1$.\\
\end{enumerate}

   \textbf{Theorem 6.16. (odd degree sine-Gordon with distinct projections)}\\
\textit{Pick any $\alpha \in \mathbb{N}^4$ satisfying $\alpha_{o}+\alpha_1\equiv 0$}(mod 2). \textit{Then, there exists a $1$-dimensional family of hyperelliptic covers $\pi: (\Gamma, p) \rightarrow (X, q)$, equipped with a pair of distinct Weierstrass points $\{p_o ,p_1\} \in \Gamma$, such that:}\\
 \begin{enumerate}
 \item $\pi(p_j)=\omega_j$,  \textit{for $j=0,1$ and} $\quad(\pi,p_o,p_1)\,\,$ \textit{solves the} \textbf{sine-Gordon case};\\
 \item $\Gamma$ \textit{has arithmetic genus} $\quad g:= \alpha^{(1)}+1\quad$ \textit{and} $\quad deg(\pi)=\alpha^{(2)}+\alpha_2+\alpha_3 +1$.\\
\end{enumerate}
   
 \textbf{Theorem 6.17. (even degree sine-Gordon with same projection)}\\
 \textit{Fix $j_o \in\{1,2,3\}$ and pick any $\alpha \in \mathbb{N}^4$ satisfying $\alpha_{j_o}+1\equiv \alpha_i$}(mod 2)\textit{ for any $i\neq j_o$. Then, there exists a $1$-dimensional family of hyperelliptic covers $\pi: (\Gamma, p) \rightarrow (X, q)$, equipped with a pair of distinct Weierstrass points $\{p_o ,p_1\} \in \Gamma$, such that:}\\
 \begin{enumerate}
 \item $\pi(p_o)=\pi(p_1)=\omega_o\quad$ \textit{and} $\quad(\pi,p_o,p_1)\,\,$ \textit{solves the} \textbf{sine-Gordon case};\\
 \item $\Gamma$ \textit{has arithmetic genus} $\quad g:= \alpha^{(1)}\quad$ \textit{and} $\quad deg(\pi)=\alpha^{(2)}+1$.\\
\end{enumerate}

 \textbf{Theorem 6.18. (odd degree sine-Gordon with same projection)}\\
 \textit{For any $\alpha \in \mathbb{N}^4$ there exists a $1$-dimensional family of hyperelliptic covers $\pi: (\Gamma, p) \rightarrow (X, q)$, equipped with a pair of distinct Weierstrass points $\{p_o ,p_1\} \in \Gamma$, such that:}\\
 \begin{enumerate}
 \item $\pi(p_o)=\pi(p_1)=\omega_o\quad$ \textit{and} $\quad(\pi,p_o,p_1)\,\,$ \textit{solves the} \textbf{sine-Gordon case};\\
 \item $\Gamma$ \textit{has arithmetic genus} $\quad g:= \alpha^{(1)}+2\quad$ \textit{and} $\quad deg(\pi)=\alpha^{(2)}+\alpha^{(1)}+3$.\\
\end{enumerate}

\noindent Universit\'e Lille Nord de France F 59000, FRANCE\\
UArtois Laboratoire de Math\'ematique de Lens EA2462,\\
F\'ederation CNRS Nord-Pas-de-Calais FR 2956\\
Facult\'e des Sciences Jean Perrin\\
Rue Jean Souvraz, S.P. 18,\\
F, 62300 LENS FRANCE\\

\noindent  Investigador PEDECIBA\\
Centro de Matem\'atica \\
Universidad de la Rep\'ublica \\
Montevideo URUGUAY


\begin{thebibliography}{99}
 
\bibitem{A-McK-M} Airault H., McKean H.P. \& Moser J., \textit{Rational and elliptic solutions of the Korteweg-deVries equation and a related many body problem}, Comm. Pure Appl. Math., \textbf{30}(1977), pp.95-148.

\bibitem{A-K-V} Akhmetshin A.A., Krichever I.M. \& Vol'vovskii Y.S., \textit{Elliptic families of solutions of the Kadomtsev-Petviashvili equation, and the field analog of the elliptic Calogero-Moser system}, Funct. Anal. and its Appl., \textbf{36} No.4 (2002), pp.253-266.

\bibitem{B-B-E} Babich M.V., Bobenko A.I. \& Enol'skii V.Z. \textit{Solutions of nonlinear equations integrable in Jacobi theta functions by the method of the inverse problem, and symmetries of algebraic curves}, Math. USSR Izvestiya vol. 26 (1986) No.2, pp.479-496.

\bibitem{B-B-E-M} Belokolos E.D., Bobenko A.I., Enol'skii V.Z. \& Matveev V. B., \textit{Algebraic-geometric principles of superposition of finite-zone solutions of the integrable nonlinear equations}, Uspekhi Math. Nauk \textbf{41} No.2 248 (1986).

\bibitem{B} Bobenko A.I., \textit{Periodic finite-zone solutions of the sine-Gordon equation}, Funct. Anal. Appl.\textbf{18} (1984).

\bibitem{C} Calogero F., \textit{Solution of the one dimensional $n$-body problems with quadratic and/or inversely quadratic pair potentials}, J. Math. Phys. \textbf{12} (1971), pp.419-436.

\bibitem{D-M-N} Dubrovin B.A., Matveev V.B. \& Novikov S.P., \textit{Nonlinear equations of KdV type, finite-band linear operators and Abelian varieties}, Uspekhi Math. Nauk \textbf{31} No.1 187 (1976). 

\bibitem{D-N} Dubrovin B.A. \& Natanzon S.M., \textit{Real two-zone solutions of the sine-Gordon equation}, Funk. Anal. Pril. \textbf{16} (1982) No.1,

\bibitem{D-N} Dubrovin B.A. \& Novikov S.P., \textit{A periodicity problem for the Korteweg-deVries-Sturm-Liouville equations. Their connection with algebraic geometry}, Doklady Acad. Nauk SSSR, \textbf{219} No.3 (1974), pp.19-22.

\bibitem{F} Fl\'edrich P.,
\textit{Paires 3-tangentielles hyperelliptiques et solutions
doublement p\'eriodiques en t de l'\'equation de Korteweg-de Vries},
Th\`ese Universit\'e d'Artois (P\^ole de Lens), D\'ec. 2003.

\bibitem {F-T} Fl\'edrich P. \& Treibich A.,
\textit{Hyperelliptic osculating covers and \textit{KdV}
solutions periodic in $t$}, I.M.R.N. \textbf{5} (2006), pp.1-17, Article ID 73476.

\bibitem{H} Hartshorne R.,
\textit{Algebraic Geometry}, Grad.Texts in Math.\textbf{52}, Springer-Verlag (1977).

\bibitem{Hi} Hirota R., \textit{Direct methods of finding exact solutions of non-linear evolution equations}, in Backlund Transformations, Springer Lecture Notes \textbf{515} (1976).

\bibitem{I-M} Its A.R. \& Matveev V.,\textit{ Hill's operator with a finite number of lacunae and
multisoliton solutions of the Korteweg-de Vries equation}, Teor. Mat. Fiz \textbf{23}
(1975), pp.51-67.

\bibitem{K-M} Kashiwara M. \& Miwa T., \textit{The} $\tau$\textit{-function of the Kadomtsev-Petviashvili equation}, Proc. Japan Acad. \textbf{57} (1981), pp.342-347.

\bibitem{K-K} Kozel V.A. \& Kotlyarov V.P., \textit{Almost periodic solutions of the equation} $u_{tt}\,$-$\,u_{xx}+sin\,u =0$, Dokl. Akad. Nauk  Ukrain SSR Ser. A \textbf{1976}, pp.878-881.

\bibitem{K.1} Krichever I.M.,
\textit{ Integration of non-linear equations by the methods of
Algebraic Geometry},
 Funct.Anal. \& Appl., \textbf{11} (1977), pp.15-31.


\bibitem{K.2} Krichever I.M.,
\textit{ Elliptic solutions of the \textit{KP} equation and
integrable systems of particles}, Funct. Anal., \textbf{14}, \textbf{4 }(1980),
pp.45-54.

\bibitem{Kru} Kruskal M.D.,\textit{ The K-deV equation and related evolution equations}, Lectures in Appl. Math. \textbf{15} Amer. Math. Soc., Providence (1974), pp.61-83.

\bibitem{Lamb} Lamb G.L., \textit{Elements of Soliton Theory}, Wiley, New York, 1980.

\bibitem{Lax} Lax P. D., \textit{Periodic solutions of the KdV equation}, Comm. Pure Appl. Math. \textbf{28} (1975), pp.141-188.

\bibitem{McK-VM} McKean H.P. \& Van Moerbeke P., \textit{The spectrum of Hill's equation}, Inv. Math. \textbf{30} (1975), pp.217-274.

\bibitem{M} Mumford D., \textit{An algebro-geometric construction of commuting operators and of solutions to the Toda lattice equation, Korteweg-de Vries equation and related nonlinear equations}, Proceedings of the International Symposium on Algebraic Geometry (Kyoto Univ., Kyoto, 1977), pp.115-153.
 
\bibitem{P} Previato E., \textit{Hyperelliptic quasi-periodic and soltion solutions fo the nonlinear Schr\"odinger equation}, Duke Math.J. Vol.52, No.2 (June 1985), pp.329-377.

\bibitem{Sato} Sato M. \& Sato Y., \textit{Soliton equations as dynamical systems on an infinite dimensional Grassmann manifold},Proc. US-Japan seminar, Non-linear partial differential equations in applied science. H.Fujita, P.D.Lax and G.Strang, ed. Konokuniya/North-Holland, 1982, pp.259-271.

\bibitem{S-W} Segal G.\& Wilson G.,\textit{ Loop groups and equations of KdV type},
Publ. Math. IHES, \textbf{61}(1985), pp.5-65.

\bibitem{S.1} Smirnov A.O., \textit{Elliptic solutions of the Korteweg-deVries equation}, Mat. Zametki \textbf{45} No. 6 (1989) pp.66-73.

\bibitem{S.2} Smirnov A.O., \textit{Finite-gap solutions of the abelian Toda Lattice of genus 4 and 5 in elliptic functions}, Teor. Mat. Fiz. \textbf{78} (1989), pp.11-21.

\bibitem{S.3} Smirnov A.O., \textit{Solutions of the KdV equation, elliptic in t},
Teor. Mat. Fiz. \textbf{100} No. 2 (1994).

\bibitem{S.4} Smirnov A.O., \textit{Elliptic solution of the nonlinear Schr\"odinger equation and modified Korteweg-de Vries equation}Russian Acad. Sci. Sb. Math. Tom 185 (1994), No.8.

\bibitem{Toda} Toda M., \textit{Waves in non-linear lattices}, Prog. Theor. Phys. Suppl., \textbf{45}, (1970), pp.174-200.

\bibitem{T.1} Treibich A.  \textit{Tangential polynomials and elliptic solitons}, Duke Math. J., \textbf{59} No.3 (1989), pp.611-627.

\bibitem{T.2} Treibich A.  \textit{Matrix elliptic solitons}, Duke Math.J., \textbf{90} No.3
(1997), pp.523-547.

\bibitem{T.3} Treibich A.  \emph{Rev\^etements hyperelliptiques d-osculateurs et solitons elliptiques de la hi\'erarchie KdV} ; C.R. Acad. Sci. Paris, S\'erie I, t. 345 (2007) 213-218. 

\bibitem{T-V.1} Treibich A.  \& Verdier J.-L., \textit{Solitons Elliptiques}; Progress in
Math.,
\textbf{88}. (app.by J.Oesterl\'e); The Grothendieck Festschrift. Ed. Birkh\"auser
(1990), pp.437-479.

\bibitem{T-V.2} Treibich A. \& Verdier J.-L., \textit{Rev\^etements exceptionnels et sommes de 4
nombres triangulaires}, Duke Math. J., vol. \textbf{68} No.2, pp.217-236.

\bibitem{Z-S} Zakharov V. E. \& Shabat A. B. \textit{A scheme for integrating the non-linear equations of mathematical physics by the methods of the inverse scattering problem}, Funct. Anal.Appl. \textbf{8} (1974), pp.226-235. \\

 \end{thebibliography}
\end{document}